\documentclass[12pt]{article}  
\usepackage{graphicx,amsmath,amsthm,amssymb}  
\numberwithin{equation}{section}
\def\RR{\mathbb R}
\def\ZZ{\mathbb Z}

\def\QQ{\mathbb Q}
\def\KK{\mathbb K}

\def\FF{\mathbb F}
\def\GG{\mathbb G}
\def\AA{\mathbb A}
\def\WW{\mathbb W}
\def\VV{\mathbb V}
\def\UU{\mathbb U}
\def\SS{\mathbb S}
\def\cK{\mathcal K}
\def\cR{\mathcal R}

\def\cA{\mathcal A}
\def\fL{\mathfrak L}
\def\cB{\mathcal B}

\def\ff{\mathfrak f}

\def\cF{\mathcal F}
\def\FF{\mathbb F}

\def\fO{\mathcal O}
\def\cD{\mathcal D}
\def\fD{\mathfrak D}
\def\cN{\mathcal N}

\def\cT{\mathcal T}
\def\fT{\mathfrak T}

\def\fU{\mathfrak U}

\def\bx{{\bf x}}
\def\by{{\bf y}}

\def\bm{{\bf m}}
\def\bu{{\bf u}}
\def\bd{{\bf d}}
\def\bN{{\bf N}}
\def\vol{{\rm vol}}
\def\bgamma{\boldsymbol{\gamma}}
\def\btheta{\boldsymbol{\theta}}
\def\bet{\boldsymbol{\eta}}
\def\bkappa{\boldsymbol{\kappa}}
\def\b0{{\bf 0}}

\def\qed{ \ \vrule width.2cm height.2cm depth0cm\smallskip}
\newcommand{\fc}{\frac}
\title{On the Gaussian limiting  distribution of lattice points in 
  a parallelepiped}
\author{Mordechay B. Levin}
\date{}
\begin{document}
\baselineskip   15pt
\hsize=15.9 true cm
\hoffset=-1.0 true cm  
\maketitle{\bf \center }
\begin{abstract}
Let $ \Gamma \subset \RR^s $ be a lattice   obtained from a
 module in a totally real algebraic number field.  Let $\cR(\btheta, \bN)$ be an error term in the lattice point problem for the parallelepiped
 $[-\theta_1 N_1,\theta_1 N_1] \times \cdots \times [-\theta_s N_s,\theta_s N_{s}]$. In this paper, we prove that $\cR(\btheta, \bN)/\sigma(\cR,\bN) $ have Gaussian limiting  distribution as $N  \to  \infty$, where
  $\btheta=(\theta_1,...,\theta_s)$ is a uniformly distributed random variable in $[0,1]^s$, $N=N_1 \cdots N_s$ and $\sigma(\cR,\bN) \asymp (\log N)^{(s-1)/2}$.
  We obtain also a similar result for the low discrepancy sequence corresponding to $\Gamma$.
\end{abstract}
Key words: lattice points problem, low discrepancy sequences, totally real algebraic number field, central limit theorem.\\
2000  Mathematics Subject Classification. Primary 11P21, 11K38, 11R80. Secondary 60F05.
%\baselineskip 15pt
%\maketitle{\bf \center }
%
%
\section{Introduction.}
\setcounter{equation}{0}
{\bf 1.1. Preliminaries.} 
In 1992, J. Beck (see [Be1]-[Be3]) discovered a very surprising phenomena of randomness of the sequence 
$\{ n  \sqrt{2}  \}_{n \ge 1}$ and the lattice $ \{ (n,n\sqrt{2}  +m) | (n, m) \in \ZZ^{2} \}$ :
\begin{equation}\nonumber
\text{vol}\left\{(x,y,z)\in [0,1)^3: \;\fc{\sum\limits_{n=0}^{[xN]}(\chi([0,y),\{ n \sqrt{2} +z\})-y)}{c_1\sqrt{\log N}}<t\right\}  \to \Phi(t)=
  \fc{1}{\sqrt{2\pi}}\int_{-\infty}^t   e^{-u^2/2}du  
\end{equation}
as $N  \to   \infty,$ where $\chi([0,y),v)$ is the indicator function of $[0,y)$, $c_1>0$ and $\{v\}$
 is the fractional part of $v$.
 
 According to [Be2, p.41],  the generalizations of this results to the multidimensional case  for a  Kronecker's lattice
 $ \{ (n,n\alpha_1 +m_1,...,n\alpha_{s-1} +   m_{s-1})
 \; | \;(n, m_1,...,m_{s-1}) \in \ZZ^{s} \}$ 
%\newpage 
\noindent
is very    difficult because of the problems connected to  Littlewood's conjecture:
\begin{equation}\nonumber
  \varliminf_{n \to \infty}   \;\;  n \ll n\alpha \gg  \ll n\beta \gg =0
\end{equation}
for all reals $\alpha, \beta$, where $\ll x\gg=\min(\{x\},1-\{x\})$. 

In this paper,  in order to avoid these  problems,  we consider a lattice $\Gamma$ obtained from a  module in a totally real algebraic number field. We prove the Central Limit Theorem for the number of points in a parallelepiped.
 We obtain also a similar result for  low discrepancy sequences corresponding to $\Gamma$ (see [Le2]).
 Results of this paper were announced in [Le1], [Le2]. 
  For  related questions and 
 generalizations, see  [Le3].
 In a forthcoming paper, we will  generalize results from [Be2] to the cases of $s$-dimensional  Halton's sequences (for $1$-dimensional case see [LeMe]), $(t,s)$-sequences, and admissible lattices (see the definition below). \\

{\bf 1.2. Lattice points.}

 Let $ \fO \subset \RR^s $ be a compact region, 
 vol $ \fO $ the volume of $\fO$, $t\fO $ the dilatation
  of $ \fO $ by a factor $ t>0$, and let $t\fO+\bx$ be the translation of 
 $t\fO $ by  a vector $ \bx \in \RR^s$.  Let $ \Gamma \subset \RR^s $ be a lattice, i.e., a
 discrete subgroup of $\RR^s $
 with a compact fundamental set $ \RR^s/\Gamma$, 
$\det \Gamma$=vol$(\RR^s/\Gamma)$.
Let
\begin{equation}\label{1.0}
  \cN (\fO ,\Gamma)={\rm card}(\fO \cap \Gamma )=\sum_{\bgamma\in
  \Gamma}\chi(\fO ,\bgamma) 
\end{equation}
be the number of points of the lattice $ \Gamma $ lying inside the
region $\fO $, where we denote by $ \chi (\fO ,\bgamma)$, $\bgamma \in
\RR^s $, the indicator function of $ \fO $. 
We define the error $\cR(\fO+\bx,\Gamma)$  by setting
\begin{equation}\label{1.2}
 \cN (\fO+\bx,\Gamma)= {\rm vol} \fO  (\det \Gamma)^{-1} \;+\;\cR(\fO+\bx,\Gamma).
\end{equation}
Let  ${\rm Nm} ( \bx)= x_1 x_2 \ldots x_s$  for $\bx=(x_1, \dots, x_s) $.
 The lattice $ \Gamma \subset \RR^s $ is {\it admissible} if
\begin{equation} \nonumber
  {\rm Nm} \;\Gamma =\inf_{\bgamma \in \Gamma \setminus \{0\}} | 
{\rm Nm} (\bgamma)| >0.
\end{equation}

Let $\cK$ be a totally real algebraic number field of degree $s \ge 2$, and let 
 $\sigma$ be the canonical embedding of $\cK$ in the Euclidean space
 $\RR^s$, $ \sigma : K\ni \xi \rightarrow \sigma (\xi) =
 (\sigma_1(\xi), \ldots, \sigma_s (\xi)) \in \RR^s $, where $
 \{\sigma_j \}_{j=1}^s $ are $s$ distinct embeddings of $\cK$ in the field $
 \RR$ of real numbers. Let  $N_{\cK/\QQ}(\xi)$ be the norm of  $\xi \in \cK$. By [BS, p. 404]
\begin{equation} \label{1.1b}
    N_{\cK/\QQ}(\xi)  = \sigma_1 (\xi) \cdots \sigma_s (\xi), \quad {\rm and}  \quad    |N_{\cK/\QQ}(\alpha)| \geq 1
\end{equation}
for all algebraic integers  $\alpha \in \cK \setminus \{ 0 \}$. We see that
$ |{\rm Nm} (\sigma (\xi))|=  | N_{\cK/\QQ}(\xi)|$. Let $M$  be a  full $\ZZ$ module  in $ \cK$, and let $\Gamma_M$ be the lattice 
 corresponding to $M$ under the embedding $\sigma$.
 It is known that the set $M^{\bot}$ of all
  $\beta \in \cK$, for which $ {\rm Tr}_{\cK/  \QQ}(\alpha \beta)  \in\ZZ$ for all $\alpha \in M$, is also a full $\ZZ$
 module ({\it the dual of the module $M$}) of the field $K$ (see [BS], p. 94).
 Recall that  the dual lattice  $ \Gamma_M^\bot $ consists of all vectors  $
\bgamma^\bot \in \RR^{s}$ such that the inner product $\langle \bgamma^\bot,\bgamma \rangle$
belongs to $\ZZ $ for each $\bgamma \in \Gamma $.  Hence $\Gamma_{M^{\bot}} =\Gamma_M^{\bot}$.
Let $(C_{M})^{-1} >0$ be an integer such that $(C_{M})^{-1} \gamma$ are algebraic integers for all $ \gamma \in M \cup  M^{\bot}$. Hence
\begin{equation} \label{1.01}
  \min( {\rm Nm}\;\Gamma_M, {\rm Nm}\;\Gamma_M^{\bot})   \geq C_{M}^s.
\end{equation} 
Therefore $\Gamma_{M} $ and $\Gamma_{M^\bot} $ are admissible lattices.
In the following we will use the notation $\Gamma=\Gamma_M$.

{\bf 1.3. Low discrepancy sequences.}

 Let $((\beta_{n,N})_{n=0}^{N-1})$ be an  $N$-point set in an $s$-dimensional unit cube $[0,1)^s$,
$\fO=[0,y_1) \times \cdots \times [0,y_s) \subseteq [0,1)^s $,
\begin{equation}\label{N1}
\Delta(\fO, (\beta_{n,N})_{n=0}^{N-1}  )= \#\{0 \leq n
<N  \;|\; \beta_{n,N}\in \fO\}-Ny_1 \ldots y_s.
\end{equation}
%\newpage
We define the $\emph{L}^\infty$ and $\emph{L}^2$ {\it discrepancy} of a 
$N$-point set $(\beta_{n,N})_{n=0}^{N-1}$ as
\begin{equation} \nonumber
  \emph{D} ((\beta_{n,N})_{n=0}^{N-1}) = 
    \sup_{ 0<y_1, \ldots , y_s \leq 1} \; \Big| \frac{1}{  N}
  \Delta(\fO,(\beta_{n,N})_{n=0}^{N-1}) \Big|,
\end{equation}
\begin{equation}\nonumber
   \emph{D}_2((\beta_{n,N})_{n=0}^{N-1})=\Big( \int_{
[0,1]^s}\Big|\frac{1}{  N}\Delta(\fO,(\beta_{n,N})_{n=0}^{N-1})\Big|^2
dy_1 \ldots dy_s \Big)^{1/2}.
\end{equation}
In 1954, Roth proved that there exists a constant $ \dot{C}>0 $, such
that
\begin{equation} \label{1.6}
N\emph{D}_2((\beta_{n,N})_{n=0}^{N-1})>\dot{C}(\ln N)^{\frac{s-1}{ 2}}, \;
 \quad    
\end{equation}
for all $N$-point sets $(\beta_{n,N})_{n=0}^{N-1}$ .

 \texttt{  Definition 1.} {\it A sequence $(\beta_n)_{n\geq 0}$ is of {\; \rm
low discrepancy} (abbreviated l.d.s.) if  $\emph{D}
((\beta_n)_{n=0}^{N-1})=O(N^{-1}(\ln N)^s) $ for $ N \rightarrow \infty $.
A sequence  of point sets $((\beta_{n,N})_{n=0}^{N-1})_{N=1}^{\infty}$ is of
 {\; \rm low discrepancy} (abbreviated
l.d.p.s.) if $ \emph{D}((\beta_{n,N})_{n=0}^{N-1})=O(N^{-1}(\ln
N)^{s-1}) $ for $ N \rightarrow \infty $.}

 For examples of l.d.s.  and  l.d.p.s. see  [BC] and [DrTi].
In [Fr], Frolov constructed a low discrepancy point set from a module in a totally real algebraic number field (see also [By], [Skr]).
 Using this approach, we proposed in [Le2] the following construction of l.d.s. :

  According to (\ref{1.01}) $|{\rm Nm} (\bgamma^{(1)} -\bgamma^{(2)})| \geq C_M^s$ for  different points $\bgamma^{(1)},\bgamma^{(2)} \in \Gamma$. Hence, there are no two different points $\bgamma^{(1)},\bgamma^{(2)} \in  \Gamma$ with 
   $\gamma_s^{(1)} = \gamma_s^{(2)}$. We have that the set $W_{\bx}=((\bx,0)+\Gamma) \cap [0,1)^{s-1} \times 
   (-\infty, \infty) $  with $\bx \in [0,1)^{s-1}  $ can be enumerated by a sequence
$ (z_{\bx,k}, z_{s}(\bx,k))_{k=-\infty}^{+\infty}$ in the following way:
\begin{equation} \label{N2}
   z_{\bx,0}= \bx,  \; z_{s}(\bx,0) =0, \;\; z_{\bx,k} \in [0,1)^{s-1}   \quad {\rm and} \quad
z_{s}(\bx,k)< z_{s}(\bx,k+1)\in \RR ,
\end{equation}
for    $k \in \ZZ$.  We see that  there exists a unique  $ (w,y_{s}) \in W_{\bx}  $ with
\begin{equation} \nonumber
y_{s} = \min \{v>0 \;|\; \exists w \in [0,1)^{s-1}, \;\; {\rm such \; that} \;\; (w,v) \in
 W_{\bx} \}.
\end{equation}
Let  $\cT (\bx)=w $. In [Le2], we proved that $\cT (\bx)$ is the ergodic transformation with respect to the 
Lebesgue measure of  $[0,1)^{s-1} $,
 $\cT^k (\bx) =z_k(\bx)$ for $k \in \ZZ$, 
   and 
  $(\cT^k (\bx))_{k \geq L}$ is the l.d.s. for all $\bx$ and all $L$.
 
We note that throughout  the paper $O$-constants do not depend on $\bx,\btheta$ and $\bN$. 
\\

{\bf 1.4.  Statement of the results.}\\

 Let $\KK_s =[-1/2,1/2)^s$, $\bN =(N_1,...,N_s)$, $N_i>0, \;(i=1,...,s)$, $N =N_1N_2...N_s $,  
$(x_1,...,x_s) \cdot (y_1,...,y_s) =(x_1y_1,...,x_s y_s)$, \; $(x_1,...,x_s) \cdot \fO =\{ (x_1,...,x_s) \cdot (y_1,...,y_s) \; : \;
 (y_1,...,y_s) \in \fO \} $ and $n = [\log_2 N]+1$.\\

{\bf Theorem 1.} {\it With the notations as above, there exist $w_2>w_1>0$ such that
\begin{equation}\nonumber
   \sup_{t, \bx} \left| \rm{vol} \left\{ \boldsymbol{\theta} \in [0,1)^s : \;\fc{ \cR( \boldsymbol{\theta} \cdot \bN \cdot \KK_s+\bx, \Gamma) }
{  w(\bN,\bx) n^{(s-1)/2}}<t \right\}  - \Phi(t) \right| =O(n^{-1/15})  
\end{equation}
 as $N  \to   \infty$, with $ w(\bN,\bx)\in [w_1, w_2]$  for all $\bx \in [0,1)^{s-1}$.}\\

{\bf Theorem 2.} {\it Let $\fO =[0,y_1) \times \cdots \times [0,y_{s-1})$. Then  there exist $w_2>w_1>0$ such that
\begin{equation}\nonumber
   \sup_t \left| \rm{vol} \left\{\by \in [0,1)^{s}, \bx \in [0,1)^{s-1}  : \;\fc{\Delta \big( \fO , (\cT^k(\bx))_{k=0}^{[y_s N]-1} \big) }
{  v(N,\bx)n^{(s-1)/2}}<t  \right\} -  \Phi(t)\right|  =O(n^{-1/15}) 
\end{equation}
 as $N  \to   \infty$, with $ v(N,\bx)\in [w_1, w_2]$.}\\
 
{\bf Remark.}  Let $K(r_1,r_2)$ be an 
algebraic number field with signature $(r_1,r_2)$, $r_1+2r_2=s$, $\Gamma=\Gamma(M,r_1,r_2) \subset \RR^s$
 a lattices obtained from a module $M$ in  $K(r_1,r_2)$, $\bN =( N_1^{'},...,N_{r_1}^{'}, N_1,...,N_{r_2} )  \in \ZZ_{+}^{r_1+r_2}$, 
 $\bgamma =( \gamma_1^{'},...,\gamma_{r_1}^{'}, \gamma_1,...,\gamma_{r_2} )  \in \RR^{s}$ $(\gamma_i^{'} \in \RR,\; \gamma_j \in \RR^2, i=1,..,r_1, j=1,...,r_2  )$,
 $\by =( y_1^{'},...,y_{r_1}^{'} , y_1,...,y_{r_2})$,   $ V=\RR^s/\Gamma$,
 $(\by,\bx)$ a uniformly distributed random variable in $[0,1]^{r_1+r_2} \times V$, ${\bf 1}_{G}$ the indicator function of the domain $G$,
 $$
      G(\bN) = \prod_{i=1}^{r_1} [-N_iy_i, N_iy_i ]  \prod_{j=1}^{r_2} \{ z \in \RR^2  \;  \: \; |z| \leq  N_j  y_j  \},
 $$
and let
$$
   \xi_1(\bN) =  \sum_{\bgamma \in \Gamma+\bx}    {\bf 1}_{G(\bN)} (\bgamma) , \quad  \xi_2(\bN) =  \sum_{\bgamma \in \Gamma+\bx}    {\bf 1}_{G(\bN)} (\bgamma) \prod_{j=1}^{r_2} \sqrt{N_j^2  y_j^2 - \gamma_j^2}
$$
 In a forthcoming paper, we will prove CLT for  the multisequence $\xi_i(\bN)$, where  $i=1$ if $r_2 \geq 2$ and  $i=2$ if $r_2=1, r_1  \geq 1$.
The case $r_2=1, r_1=0$ was investigated earlier by  Hughes and Rudnick [HuRu]. \\

\section{ Proofs of theorems}
\setcounter{equation}{0}
\renewcommand{\theequation}{2.1.\arabic{equation}}
{\bf 2.1. Algebraic units.} 

   Let $\fD_M$ be the ring of coefficients of the full module $M$,  $\fU_M$ be the group of units of $\fD_M$, $M_1=M$, $M_2=M^{\bot}$,
and let   $\eta_{k,1},..., \eta_{k,s-1}$ be the set of fundamental units of $\fU_{M_k}\; (k=1,2)$. According to the Dirichlet's
 theorem (see e.g., [BS, p. 112]), every unit $\eta \in \fU_{M_k}$ has a unique representation in the form
\begin{equation}\label{2.1}
   \eta = (-1)^a\eta_{k,1}^{a_1} ...\eta_{k,s-1}^{a_{s-1}},  \quad \quad k=1,2,
\end{equation} 
where $a_1,...a_{s-1}$ are rational integers and $a \in \{0,1\}$.   We  will denote  $\sigma(\fU_{M_k})$  by the same
 symbol $\fU_{M_k}$.\\

{\bf Lemma 1.} {\it Let $y_1,...,y_s >0$ be reals, $\by=(y_1,...,y_s)$, $y= {\rm Nm} ( \by) = y_1y_2...y_s$. 
Then there exists $\eta_k(\boldsymbol{y}) \in \fU_{M_k}$ with}
\begin{equation}\label{2.2}
   y_i y^{-1/s}|\sigma_i(\eta_k(\boldsymbol{y}) )| \in [1/c_0, c_0],  
\end{equation} 
where $ i=1,...,s, \; k=1,2,$ and
\begin{equation}\label{2.2a}
                c_0 = \exp\Big( \sum_{k =1,2}\sum_{1 \leq i,j <s}  |\ln |\sigma_i(\eta_{k,j})||\Big) >1.
\end{equation} \\

 {\bf Proof.} We fix $k \in \{1,2 \}$. The matrix $ (\ln |\sigma_i(\eta_{k,j})|)_{1 \leq i,j <s}$ is non singular [BS, pp. 104,115]. Hence, there exist reals $b_1,...,b_{s-1}$ with
\begin{equation} \nonumber
     \sum_{1 \leq j <s} b_j \ln |\sigma_i(\eta_{k,j})|  =1/s \ln y -  \ln y_i, \quad \quad  i=1,...,s-1.   
\end{equation}  
Taking $a_j=[b_j], \; j=1,...,s-1$ and $ \eta_k(\boldsymbol{y}) = \eta_{k,1}^{a_1} ...\eta_{k,s-1}^{a_{s-1}} $, we obtain for $i \in [1,s-1]$
\begin{equation}\label{2.4}
  -\sum_{1 \leq j <s} | \ln |\sigma_i(\eta_{k,j}) || \leq
         \ln (y_i y^{-1/s}|\sigma_i(\eta_k(\boldsymbol{y}) )|) 
         \leq  \sum_{1 \leq j <s} | \ln |\sigma_i(\eta_{k,j}) || .
\end{equation} 
Hence 
\begin{equation}\nonumber
 -\ln c_0   \leq  \ln (y_i y^{-1/s} |\sigma_i(\eta_k(\by) )|) \leq \ln c_0, \quad i=1,...,s-1, \;\; k=1,2 .
\end{equation} 
Bearing in mind that $|{\rm Nm} (\eta_k( \by))| =1$ and  $y = y_1y_2...y_s $, we get from   (\ref{2.4}) and (\ref{2.2a})
\begin{equation}\nonumber
         \ln (y_s y^{-1/s}|\sigma_s(\eta_k(\boldsymbol{y}) )|) 
          = -\sum_{1 \leq i <s}   \ln (y_i y^{-1/s}|\sigma_i(\eta_k(\boldsymbol{y}) )|) \in [ -\ln c_0  ,\ln c_0 ] .
\end{equation}  
Therefore, the assertion (\ref{2.2}) is true for $i \in [1,s], \; k=1,2$, and  Lemma 1 is proved. \qed \\

  We apply this lemma to the vector $\by=\bN= (N_1,...,N_s)$.   
   Let $N_i^{'} =  N_i |\sigma_i(\eta_1(\bN))|$,\\  $i=1,...,s$
  and let $\sigma(\eta_1(\bN)) = ( \sigma_1(\eta_1(\bN)),...,\sigma_s(\eta_1(\bN)))$.
We see that 
\begin{equation}\nonumber
         \bgamma \in \Gamma_M     \cap  (\btheta \cdot \bN \cdot \KK^s +\bx)  \Leftrightarrow \bgamma \cdot \sigma(\eta_1(\bN))  \in \Gamma_M   
           \cap  (\btheta \cdot \bN^{'} \cdot \KK^s +\bx \cdot \sigma(\eta_1(\bN) )).
\end{equation}  
Hence
\begin{equation} \nonumber
   \cN(\btheta \cdot \bN \cdot \KK^s +\bx  ,\Gamma_M  )  = \cN(\btheta \cdot \bN' \cdot \KK^s +\bx \cdot \sigma(\eta_1(\bN) ),\Gamma_M  ).
\end{equation}  
By (\ref{1.2}), we have 
\begin{equation}\nonumber
   \cR(\btheta \cdot \bN \cdot \KK^s +\bx  ,\Gamma_M  )  = \cR(\btheta \cdot \bN' \cdot \KK^s +\bx \cdot \sigma(\eta_1(\bN)) ,\Gamma_M  ).
\end{equation}  
Therefore, without loss of generality, we can assume that
\begin{equation}\label{2.81}
    N_i N^{-1/s} \in [1/c_0, c_0], \quad  i=1,...,s.  
\end{equation}  
Now, let $n=[\log_2 N] +1$,
\begin{equation} \nonumber
       \FF_n^{'} = \{\bgamma \in \Gamma^{\bot}  \setminus \{0 \} \; : \; 
               |\gamma_i | |{\rm Nm} ( \bgamma)|^{-1/s}  \in [1/c_0,c_0],  \; i=1,...,s  ,  \quad |{\rm Nm} ( \bgamma)|  \leq n^{1/2} \}         
\end{equation} 
and
\begin{equation}\label{2.82}
     \FF_n = \bigcup_{\bgamma \in  \FF_n^{'}} \{\bgamma^{'} \in  \FF_n^{'}  \; : \; \gamma^{'}_1 = 
     \max_{ \eta \in\fU_{M^\bot}, \;  \bgamma \cdot \sigma( \eta) \in  \FF_n^{'} } (\bgamma \cdot \sigma( \eta))_1 \}.
\end{equation}   
By  (\ref{1.01}),  we get that
%  Nigde ne nujno ?
\begin{equation}\label{2.82a}
   {\rm if} \quad \bgamma^{(1)},\bgamma^{(2)} \in  \FF_n, \;  \bgamma^{(1)} \neq \bgamma^{(2)}   \quad  \quad {\rm then}  \quad  \quad 
   \bgamma^{(1)} \neq \bgamma^{(2)}  \cdot \sigma( \eta)     \quad \forall \eta \in\fU_{M^\bot}.
\end{equation}      
\\

{\bf Lemma 2.} {\it Let  $a, b \geq 1 $ be integers,
\begin{equation} \label{2.90}
     \GG(a,b) = \{\bgamma \in \Gamma^{\bot}  \setminus \{0 \} \; : \;  \max_{1 \leq i \leq s}|\gamma_i| \in
       ( 2^a, 2^{a+b}], \quad
              |{\rm Nm}( \bgamma )| \leq n^{1/2}    \},
\end{equation} 
\begin{equation}\label{2.91}
     \GG^{'}(a,b):= \bigcup_{\bgamma \in    \FF_n } \bigcup_{ \eta \in  \fU(\bgamma,a,b) } \bgamma \cdot \sigma(\eta),
\end{equation}   
with
\begin{equation}\label{2.9}
  \dot{\fU}(\bgamma^{(0)},a,b) = \{ \eta \in\fU_{M^\bot} \; : \;  \bgamma^{(0)}\cdot \sigma(\eta)  \in \GG(a,b) \}.
\end{equation}  
Then}
\begin{equation}\label{2.10}
   \GG(a,b) =   \GG^{'}(a,b),  \quad   \quad      \# \dot{\fU}(\bgamma^{(0)}, a,b) =O(b(a+b)^{s-2}) 
         \quad {\rm for}  \quad \bgamma^{(0)} \in    \FF_n, 
\end{equation} 
\begin{equation}  \nonumber
   \# \GG(a,b) =   O(n^{1/2}b(a+b)^{s-2}), 
\end{equation}  
\begin{equation} \label{2.10a}
  \sum_{\bgamma  \in \FF_n} \frac{1}{|{\rm Nm}( \bgamma )|} =O(\ln n)   \qquad {\rm and} \qquad
   \sum_{\bgamma  \in \FF_n} \frac{1}{{\rm Nm}^2( \bgamma )} =O(1).
\end{equation}  \\

{\bf Proof.} It is easy to see that $\GG(a,b) \supseteq \GG^{'}(a,b)$. Let $ \bgamma  \in \GG(a,b)$. By Lemma 1, there exists $ \eta \in \fU_{M^\bot}$
with $ \bgamma \cdot \sigma(\eta)  \in   \FF_n^{'}$. From (\ref{2.82}), we obtain that there exists  $ \eta_1 \in \fU_{M^\bot}$  with 
 $\bgamma^{(1)}=\bgamma \cdot \sigma(\eta \eta_1 ) \in    \FF_n$.
By  (\ref{2.91}) and  (\ref{2.9}), we get that $\bgamma =\bgamma^{(1)} \sigma((\eta \eta_1)^{-1} )\in \GG^{'}(a,b)$ and $\GG(a,b) = \GG^{'}(a,b)$. 
  
 Let $ \bm=(m_1,...,m_{s}) \in \ZZ^s$, $m_1+...+m_{s}=0$, $\bkappa=(\kappa_1,...,\kappa_s)$, $\kappa_i \in \{-1,1 \} \; (i=1,...,s)$,
 $  \nu(\mu) =s$ if $\mu \neq s$,  $  \nu(\mu) =1$ if $\mu = s$, $j \geq 0$ and
\begin{equation}\nonumber
 B(\bm,\mu, \bkappa,j)  = \prod_{1 \leq i  <  \nu(\mu)} (\kappa_i2^{m_i},\kappa_i2^{m_i+1} ] \times 
     (j\kappa_{ \nu(\mu)}2^{-m_{ \nu(\mu)}} C_M^s   ,  (j+1)\kappa_{ \nu(\mu)}2^{-m_{ \nu(\mu)}} C_M^s  ]
\end{equation}  
\begin{equation}\label{2.102}
    \times \prod_{ \nu(\mu) <i  \leq s} (\kappa_{i}2^{m_i},\kappa_{i}2^{m_i+1} ] .
\end{equation} 
It is easy to see  that
\begin{equation}\label{2.102aa}
    B(\bm_1,\mu, \bkappa_1,j_1) \cap   B(\bm_2,\mu, \bkappa_2,j_2) = \emptyset \quad {\rm for}  \quad  (\bm_1,\mu, \bkappa_1,j_1)
                  \neq (\bm_2,\mu, \bkappa_2,j_2).
\end{equation} 
Applying (\ref{1.01}),   we have for  every $\mu \in [1,s]$ that
\begin{equation} \label{2.102a}
   \Gamma^{\bot} \setminus \{ 0 \}  =   \bigcup_{\kappa_1,...,\kappa_s \in \{-1, 1\}}   \bigcup_{\bm \in \ZZ^s, \atop m_1+...+m_s=0} 
\bigcup_{j \geq 0}
      \bigcup_{ \bgamma \in B(\bm,\mu,  \bkappa,j) } \bgamma .
\end{equation}
Let
\begin{equation}\nonumber
  \bgamma^{(1)},\bgamma^{(2)} \in \Gamma^{\bot}   \cap  B(\bm,\mu, \bkappa,j) .
\end{equation}
From (\ref{2.102}),   we see that
$$
    | {\rm Nm}(\bgamma^{(1)}  - \bgamma^{(2)}   )| < C_M^s.
$$
By  (\ref{1.01}),  we obtain that   $\bgamma^{(1)}=\bgamma^{(2)}$ and
\begin{equation}\label{2.121}
                  \#\Gamma^{\bot}   \cap  B(\bm, \mu,\bkappa,j) \leq 1.
\end{equation} 
Suppose 
$$
   \eta \in \fU_{M^{\bot}}   \cap  B(\bm,\mu, \bkappa,j).
$$
Using  (\ref{2.102}), we have that
\begin{equation}\label{2.13}
  1= |{\rm Nm}( \eta)|= (j+z_1)C_M^s 2^{z_2(s-1)},  \quad \quad {\rm with}  \quad \quad z_1,z_2 \in [0,1].
\end{equation}
Hence 
\begin{equation}\nonumber
      -1 +2^{1-s} /C_M^s     \leq j \leq   1/C_M^s.
\end{equation}
Applying (\ref{2.121}), we get
\begin{equation}\label{2.141}
               \sum_{j  \geq 0}  \# \fU_{M^{\bot}}   \cap  B(m, \mu,\bkappa,j)  \leq 2+1/C_M^s.
\end{equation}
We denote $\sigma^{-1}(B(m, \mu,\bkappa,j)) $ and $ \sigma^{-1}(\dot{\fU}(\bgamma^{(0)},a,b))  $ by the same symbols  $B(m, \mu,\bkappa,j) $  and $ \dot{\fU}(\bgamma^{(0)},a,b)  $.
Now let 
\begin{equation} \nonumber 
    \ddot{\fU}_{\mu}(\bgamma^{(0)},a,b) = \{ \bgamma \in   \dot{\fU}(\bgamma^{(0)},a,b) \;\; : \;\; |\gamma_i| \leq  |\gamma_{\mu}^{(0)}|,
      \quad i=1,...,s   \}.
\end{equation}
It is easy to see that
\begin{equation} \label{2.141a}
    \dot{\fU}(\bgamma^{(0)},a,b)  = \bigcup_{\mu \in [1,s]}\ddot{\fU}_{\mu}(\bgamma^{(0)},a,b). 
\end{equation}
Let 
\begin{equation}\label{2.142}
    \eta \in  \ddot{\fU}_{\mu}(\bgamma^{(0)},a,b)  \cap  B(m, \mu,\bkappa,j).
\end{equation}
 Denote  $m_i \in \ZZ \; (i=1,...,s)$ from the following condition:  
\begin{equation}\label{2.142a}
 \log_2 | \sigma(\eta)_{i}|=  m_{i} + z_i,     \quad \quad {\rm with}  \quad \quad z_i \in [0,1).
\end{equation} 
By  (\ref{2.90}) and  (\ref{2.142}), we obtain 
\begin{equation}\nonumber
 \log_2 |(\bgamma^{(0)} \sigma(\eta))_{\mu}|= \log_2 |\gamma_{\mu}^{(0)}| + m_{\mu} + z_{\mu}  \in (a,a+b] ,
\end{equation}
and
\begin{equation}\label{2.144}
      m_{\mu} \in J_1 := (a-1 - \log_2 |\gamma_{\mu}^{(0)}|,a+b- \log_2 |\gamma_{\mu}^{(0)}|] \cap \ZZ,    
        \quad  {\rm with}  \quad \#J_1 \leq b+2.
\end{equation}
From   (\ref{2.90}), (\ref{2.142}) and (\ref{1.01}), we get
\begin{equation}\label{2.141b}
  \log_2 |\gamma_i^{(0)} \sigma(\eta)_i| =   \log_2 |\gamma_i^{(0)}|+    m_i +z_i \leq a+b 
        \quad \quad {\rm and }  \quad \quad  m_i \leq a+b -  \log_2 |\gamma_i^{(0)}|.
\end{equation}
Using (\ref{2.141b}), (\ref{2.142a}) and (\ref{2.142a}), we derive that
\begin{equation}\nonumber
    \log_2 |\gamma_{i}^{(0)} \sigma(\eta)_{i} | = \log_2 |\gamma^{(0)}_{i} |  -\sum_{j \in [1,s],j \neq i}  \log_2|\sigma(\eta)_i|
\end{equation}
\begin{equation}\nonumber
  \geq  \log_2 |\gamma^{(0)}_{i} |  -\sum_{j \in [1,s],j \neq i}  (m_j+1) \geq  \sum_{j \in [1,s]}  \log_2|\gamma^{(0)}_j|  -(s-1)(a+b+1) 
\end{equation}
\begin{equation}\nonumber
 =    -(s-1)(a+b+1) + \log_2|{\rm Nm}(\bgamma^{(0)})| \geq  -(s-1)(a+b+1) + \log_2C_M^s.
\end{equation}
By (\ref{2.141b}), we have   $     m_i \in [\log_2 |\gamma_i^{(0)} \sigma(\eta)_i| -  \log_2 |\gamma_i^{(0)}| -1, a+b -  \log_2 |\gamma^{(0)}_i|]$. 
Hence
\begin{equation}\nonumber
  m_i \in J_2 :=[-1 -(s-1)(a+b+1) +s\log_2C_M -  \log_2 |\gamma^{(0)}_i|,a+b -  \log_2 |\gamma^{(0)}_i|],
\end{equation}
with  $ \#J_2 \leq s(a+b+1) +2 +s|\log_2 C_M|$.

 We fix $ \mu \in [1,s]$ and we consider (\ref{2.102a}). For given $m_{1},...,m_{\nu(\mu) -1},m_{\nu(\mu) +1},...,m_s  $,  we take  $m_{\nu(\mu) } =- \sum_{i \in [1,s], i \neq  \nu(\mu)  } m_i$.
  By (\ref{2.102a}),   we get 
\begin{equation} \nonumber
    \# \ddot{\fU}_{\mu}(\bgamma^{(0)},a,b) \leq   \sum_{\kappa_1,...,\kappa_s \in \{-1, 1\}}   \sum_{m_{\mu}\in J_1} 
    \sum_{m_i \in J_2 \atop i \neq \mu, \nu(\mu)  }  \;\sum_{j  \geq 0}
      \#( \dot{\fU}(\bgamma^{(0)},a,b)   \cap  B(m,\mu,  \bkappa,j) ) .
\end{equation}
Bearing in mind (\ref{2.141}), (\ref{2.141a})  and (\ref{2.144}),   we obtain 
\begin{equation}  \label{2.141z}
       \# \dot{\fU}(\bgamma^{(0)},a,b)  =O( \#J_1 (\#J_2)^{s-2}) = O(b(a+b)^{s-2}). 
\end{equation}
Hence, the assertion (\ref{2.10}) is proved. 

Let $F_1 \subset \RR^s$ be a fundamental domain for the field $\cK$, and let 
$$
F_2 = \{\bgamma \in \Gamma^{\bot}  \setminus \{0 \} \; : \; 
               |\gamma_i | |{\rm Nm} ( \bgamma)|^{-1/s}  \in [1/c_0,c_0],  \quad i=1,...,s \} \qquad ({\rm see} \; (\ref{2.2}) )  .
$$
 By [BS, pp. 312, 322], the points of $F_1$ can be arranged in a sequence $\dot{\bgamma}^{(k)}$ so that $0 < |{\rm Nm}(\dot{\bgamma}^{(1)}) | \leq  |{\rm Nm}(\dot{\bgamma}^{(2)})| \leq ... $
 and $c^{(1)} k   \leq |{\rm Nm}(\dot{\bgamma}^{(k)})|  \leq c^{(2)}  k$,   $k=1,2,...$ for some $c^{(2)} > c^{(1)}>0$. 
 Therefore, the points of $F_2$ can be arranged in a sequence $\bgamma^{(k)}$ so that $0 < |{\rm Nm}(\bgamma^{(1)}) | \leq  |{\rm Nm}(\bgamma^{(2)})| \leq ... $ and
\begin{equation}  \nonumber
                  c^{(3)} k   \leq |{\rm Nm}(\bgamma^{(k)})|  \leq c^{(4)}  k,  \quad k=1,2,...
\end{equation}
for some $c^{(4)} > c^{(3)}>0$. 

Using  (\ref{2.82}), we have that
\begin{equation}  \nonumber
                  \sum_{\bgamma  \in \FF_n} 1/|{\rm Nm}( \bgamma )| =O(\ln(n)), \quad 
                  \sum_{\bgamma  \in \FF_n} 1/{\rm Nm}^2( \bgamma ) =O(1), \quad {\rm and}  \quad \#\FF_n = O(n^{1/2}).
\end{equation}
By (\ref{2.91}) and (\ref{2.141z}), we obtain 
\begin{equation} \nonumber
       \# \GG(a,b) \leq \sum_{\bgamma^{(0)} \in  \FF_n} \; \#\dot{\fU}(\bgamma^{(0)},a,b) =O(n^{1/2}b(a+b)^{s-2}) .
\end{equation} 
Hence, Lemma 2 is proved. \qed  
\\

\setcounter{equation}{0}
\renewcommand{\theequation}{2.2.\arabic{equation}}
{\bf 2.2. Diophantine inequalities.}

 We consider the following simple variant of the \texttt{ S-unit theorem} (see  [ESS,  Theorem~1.1, p. 808]):   
Let $\beta_1,...,\beta_d \in \cK$, $\beta_i \neq 0, \; i=1,..,d$, $\deg(\cK)=s$. We consider the equation
\begin{equation}\label{6.0}
      \beta_1 \eta_1 +... +\beta_d \eta_d =1
\end{equation}
with $\bet=(\eta_1,...,\eta_d) \in (\fU_{M^{\bot}})^d$.  A solution $\bet$ of (\ref{6.0})
 is called \texttt{ non-degenerate}
if  $ \sum_{i \in I}  \beta_i \eta_i  \neq 0$ for every nonempty subset $I$ of $\{1,...,d\}$.  \\

{\bf Theorem  A.} {\it The number  $A(\beta_1,...,\beta_d)$ of non-degenerate solutions $\bet \in (\fU_{M^{\bot}})^d$ of equation (\ref{6.0})
satisfies the estimate }
\begin{equation}\nonumber
     A(\beta_1,...,\beta_d) \leq \exp((6d)^{3d}s).
\end{equation}
\\

\texttt{   Linear forms in logarithms.}
Write $\Lambda$ for the linear form in logarithms,
\begin{equation}\nonumber
    \Lambda =b_1 \log \alpha_1 + ...+ b_k \log \alpha_k     ,
\end{equation}
where $b_1,...,b_k$ are integers, $|b_i| \leq B\; (i=1,...,k)$, $B \geq e$. We shall assume  that $\alpha_1,...,\alpha_k $
 are non-zero algebraic numbers
with  heights at most $A_1,...,A_k$ (all $\geq e$) respectively.  \\

 {\bf Theorem B.} [BW, Theorem 2.15, p. 42] {\it  If $\Lambda \neq 0$, then
\begin{equation}\nonumber
   | \Lambda |  >\exp( -(16kd)^{2(k+2)} \ln A_1 ... \ln A_k  \ln B) ,
\end{equation}
where $d$ denote the degrees of} $\QQ(\alpha_1,...,\alpha_k )$. 
\\

 Let
\begin{align}
    &G^{(1)} = \{\bgamma \in \Gamma^{\bot}  : |\bgamma | \leq N,  \; | {\rm Nm}(\bgamma ) |\leq n^{1/2} 
         \;\;{\rm and}  \;\;     |N_i \gamma_i| > 2^{(\ln n)^4}  \;\; \forall i \in [1,s]  \} ,           \label{3.21a}   \\   
 &G^{(2)} =\{ \bgamma  \in \Gamma^{\bot}  :   |\bgamma |>N^5 \}, \quad 
      G^{(3)} =\{\bgamma  \in \Gamma^{\bot}  :   |\bgamma | \leq N^5, \;  |{\rm Nm}(\bgamma ) |> n^{1/2}  \},  \label{3.21} \\
 &G^{(4)} =\{\bgamma \in \Gamma^{\bot} \; :  \; N <|\bgamma | \leq N^5,  \; |{\rm Nm}(\bgamma )|\leq n^{1/2}  \}
          \;{\rm and} \;   G^{(5)} = \{\bgamma \in \Gamma^{\bot}  \setminus \{0 \} \; :  \nonumber  \\ 
  & \qquad   \quad |\bgamma | \leq N,  \; | {\rm Nm}(\bgamma )  |\leq n^{1/2}   
               \quad {\rm and}  \quad   \exists  i \in [1,s] \quad {\rm with} \quad          |N_i \gamma_i| \leq 2^{(\ln n)^4}      \}. \label{3.22}    
\end{align}
It is easy to see that $G^{(i)}  \cap  G^{(j)} =\emptyset$ for $i \neq j$, and
\begin{equation}\label{3.20}
    \Gamma^{\bot}  \setminus   \{0 \}  =       G^{(1)}  \cup  G^{(2)}  \cup  G^{(3)}    \cup  G^{(4)} \cup  G^{(5)}.
\end{equation}
Let
\begin{equation}\label{3.56}
      \dot{G}_0 = \{\bgamma \in G^{(1)}   \; : \;                  \quad     \max_{1 \leq j \leq s} |\gamma_i|  \leq 2^{n^{4/9}}  \},                 
\end{equation}
\begin{equation}\nonumber
      \dot{G}_i = \{\bgamma \in G^{(1)}   \; : \; 
                 \quad     \max_{1 \leq j \leq s} |\gamma_i| \in (2^{i n^{4/9}}, 2^{(i+1) n^{4/9} - n^{2/9}} ]   \},                   
\end{equation}
and
\begin{equation}\nonumber
      \ddot{G}_i = \{\bgamma \in G^{(1)}   \; : \; 
                 \quad     \max_{1 \leq j \leq s} |\gamma_i| \in ( 2^{(i+1) n^{4/9} - n^{2/9}},2^{(i+1) n^{4/9}}  ]   \},  \quad i=1,2,...     
\end{equation}
By (\ref{3.21a}) and (\ref{3.56}),  we have that  $\dot{G}_i  \cap  \dot{G}_j =\emptyset$ for $i \neq j$, $\dot{G}_i  \cap  \ddot{G}_j =\emptyset$ and
\begin{equation}\label{3.57}
    G^{(1)} =    \dot{G}_0 \cup \bigcup_{i=1}^{[n^{5/9}]}      (  \dot{G}_i \cup \ddot{G}_i) .
\end{equation}
\\

 {\bf Lemma 3.} {\it There exist   $\dot{c},\ddot{c}>0$ such that for all $\nu \in [1,s] $ and $\kappa \in \{ -1,1 \} $ }
\begin{equation} \nonumber
    \min_{\bgamma^{(1)},\bgamma^{(2)} \in G^{(1)}, \gamma_{\nu}^{(1)} \neq \kappa \gamma_{\nu}^{(2)}}  
  N_{\nu}|\gamma_{\nu}^{(1)} -\kappa \gamma_{\nu}^{(2)}|  \geq
   N_{\nu} |\gamma_{\nu}^{(2)}| \exp( -\ddot{c}  (\ln n)^3) \geq  \dot{c} n^{20s}.
\end{equation}

{\bf Proof.} Let  $\gamma_{\nu}^{(1)}/\gamma_{\nu}^{(2)} \kappa <0$. From (\ref{3.21a}), we obtain 
\begin{equation} \label{6.100a}
    \min_{\bgamma^{(1)},\bgamma^{(2)} \in G^{(1)}, \gamma_{\nu}^{(1)} \neq \kappa \gamma_{\nu}^{(2)}}  
  N_{\nu}|\gamma_{\nu}^{(1)} -\kappa \gamma_{\nu}^{(2)}|  \geq
   N_{\nu} |\gamma_{\nu}^{(2)}|  \geq  2^{(\ln n)^4}.
\end{equation}
 Now let
$\gamma_{\nu}^{(1)}/\gamma_{\nu}^{(2)}  \kappa >0$. Taking into account that $|\exp(x) -1| \geq |x|$ for any real $x$, we get
\begin{equation} \label{6.100}
    |\gamma_{\nu}^{(1)} -\kappa \gamma_{\nu}^{(2)}| =   |\gamma_{\nu}^{(2)}(\exp(\ln(\kappa\gamma_{\nu}^{(1)} /\gamma_{\nu}^{(2)})) -1) | 
      \geq    |\gamma_{\nu}^{(2)} \ln(\kappa \gamma_{\nu}^{(1)} /\gamma_{\nu}^{(2)})|.
\end{equation}
   By (\ref{2.82}), (\ref{2.90}), (\ref{2.91}) and (\ref{3.21a}),  we have
 that there exists $(\dot{\gamma}_{\nu}^{(k)},\eta_{k} )$  such that $ \gamma_{\nu}^{(k)} =  \dot{\gamma}_{\nu}^{(k)}    \cdot \sigma_{\nu}(\eta_{k})$ where $\dot{\gamma}^{(k)} \in \FF_n$ and $\eta_k$ is a unit in 
 $\cK \; (k=1,2)$. Let $\dot{\gamma}^{(k)} =\sigma(\ff^{(k)})$ with some $\ff^{(k)} \in M^{\bot}   \; (k=1,2)$.
Using  (\ref{2.1}), (\ref{2.82}) and (\ref{3.21a}), we obtain
\begin{equation}\nonumber
         \gamma_{\nu}^{(k)} =  \sigma_{\nu}(\ff^{(k)})    (-1)^{a^{(k)}} \sigma_{\nu}(\eta_{2,1})^{a_1^{(k)} } \cdots   
     \sigma_{\nu}(\eta_{2,s-1})^{a_{s-1}^{(k)} },
\end{equation}
\begin{equation}\label{6.100f}
       |\sigma_{i}(\ff^{(k)})|  \leq c_0 n^{1/(2s)}  \quad {\rm for} \quad i=1,...,s,
\end{equation}
and
\begin{equation} \nonumber
          |a_1^{(k)} \ln(\sigma_{\nu}(\eta_{2,1})) + ... +  a_{s-1}^{(k)} \ln(\sigma_{\nu}(\eta_{2,s-1})) | \leq 
\end{equation}
\begin{equation} \nonumber
      \leq           | \ln|\gamma_{\nu}^{(k)}|| + |\ln|\sigma_{\nu}(\ff^{(k)}) ||  
       \leq   \ln N + 1/(2s) \ln n + \ln (c_0). 
\end{equation}
Bearing in mind that $\det((\ln(\sigma_{i}(\eta_{2,j})))_{1 \leq i,j \leq s-1} ) \neq 0$ (see [BS, pp. 104, 115]), we get that there exists $ \tilde{C}_{1}>0$ such that 
$|a_{i}^{(k)}| < \tilde{C}_{1} n$ for $i=1,...,s-1, k=1,2$ and $n=[\log_2 N]+1$.

Let $\kappa_1 ={\rm sign}(\gamma_{\nu}^{(1)} /\gamma_{\nu}^{(2)})$, where ${\rm sign}(x)=1$ for $x>0$ and  ${\rm sign}(x)=-1$ for $x<0$   . We see that
\begin{equation} \nonumber
 \ln|\gamma_{\nu}^{(1)} /\gamma_{\nu}^{(2)}|=\ln( \kappa_1 \gamma_{\nu}^{(1)} /\gamma_{\nu}^{(2)}) 
 = \ln( \kappa_1(-1)^{a^{(1)}}\sigma_{\nu}(\ff^{(1)})   ) - \ln( (-1)^{a^{(2)}}\sigma_{\nu}(\ff^{(2)})  )
\end{equation}
\begin{equation} \nonumber
    +(a_1^{(1)} - a_1^{(2)}) \ln(\sigma_{\nu}(\eta_{2,1})) + \cdots + (a_{s-1}^{(1)} - a_{s-1}^{(2)}) \ln(\sigma_{\nu}(\eta_{2,s-1})).
\end{equation}
Let  $ \tilde{C}_{2} \max_{i \in [1,s-1]} H(\eta_{2,i})$, where  $H(\alpha)$ is the height of $\alpha$. 
By  (\ref{1.01}), $C_{M}^{-1} \ff^{(k)}$ is an algebraic integer $(k=1,2)$. Thus 
$f(x) =x^s + f_{s-1} x^{s-1}+\cdots +f_0=(x- \sigma_{1}(C_{M}^{-1}\ff^{(k)}) ) \cdots  \\(x- \sigma_{s}(C_{M}^{-1}\ff^{(k)}) )$ is the characteristic polynomial of $C_{M}^{-1} \ff^{(k)}$.
 Hence \\$H(C_{M}^{-1} \ff^{(k)}) \leq \max_{i \in [0,s-1]} |f_i|$. 
From (\ref{6.100f}), we have that  $H(C_{M}^{-1} \ff^{(k)}) \leq (2 C_{M}^{-1} c_0  n^(1/(2s))))^s$ and   $H( \ff^{(k)}) \leq (2 C_{M}^{-2} c_0  n^(1/(2s))))^s$.\\
Applying Theorem B  with $d=s, k=s+1$, 
 $\alpha_1 =\kappa_1(-1)^{a^{(1)}}\sigma_{\nu}(\ff^{(1)}) $,  $\alpha_2 =(-1)^{a^{(2)}}\sigma_{\nu}(\ff^{(2)}) $,
  $\alpha_3 =\sigma_{\nu}(\eta_{2,1})$,..., $\alpha_{s+1} =\sigma_{\nu}(\eta_{2,s-1})$, 
 $ A_1=A_2=(2 C_{M}^{-2} c_0  n^(1/(2s)))^s, A_3=\cdots =A_{s-1}= \tilde{C}_{2}$ and $B=  2 \tilde{C}_{1}n$,   we obtain
\begin{equation} \nonumber
 |\ln( \kappa_1 \gamma_{\nu}^{(1)} /\gamma_{\nu}^{(2)}) |  \geq  \exp( -\ddot{c}_{\nu}  (\ln n)^3),
\end{equation}
with some $\ddot{c}_{\nu}  >0$.
Taking into account (\ref{6.100}) and that $N_{\nu}  |\gamma_{\nu}^{(2)}| \geq  2^{(\ln  n)^4} $, we have 
\begin{equation} \nonumber
    \min_{\bgamma^{(1)},\bgamma^{(2)} \in G^{(1)}, \gamma_{\nu}^{(1)} \neq \kappa \gamma_{\nu}^{(2)}}  
  N_{\nu}|\gamma_{\nu}^{(1)} -\kappa \gamma_{\nu}^{(2)}|  \geq
   N_{\nu} |\gamma_{\nu}^{(2)}| \exp( -\ddot{c}_{\nu}   (\ln n)^3) \geq  \dot{c}_{\nu}  n^{20s},
\end{equation}
with some $\dot{c}_{\nu}  >0$.  Now using (\ref{6.100a}), we get the assertion of Lemma 3. \qed \\

{\bf 2.3. Poisson summation formula.} 

\setcounter{equation}{0}
\renewcommand{\theequation}{2.3.\arabic{equation}}
We shall need the Poisson summation formula: 
\begin{equation}\label{10.0}
             \det \Gamma   \sum_{\bgamma\in \Gamma} f(\bgamma-X)  =   \sum_{\bgamma\in \Gamma^\bot }  \widehat{f}(\bgamma)  e(\langle\bgamma,\bx\rangle),
\end{equation}
where  
\begin{equation} \nonumber
            \widehat{f}(Y) = \int_{\RR^s}{ f(X) e(\langle \by,\bx\rangle) d\bx} 
\end{equation}
is the Fourier transform of $f(X)$, and $e(x) =exp(2\pi \sqrt{-1}x), \langle\by,\bx\rangle = y_1x_1+ \cdots +y_s x_s$. Formula (\ref{10.0}) holds for functions $f(\bx)$ with period lattice $\Gamma$ 
if one of the functions $f$ or $\widehat{f}$ is integrable  and belongs to class $C^{\infty}$ (see e.g.  [StWe, p. 251]).

Let $\bd =(d_1,...,d_s)$, $d_i \geq 0$ $(i=1,...,s)$, $\fO_{\bd} =[-d_1/2, d_1/2] \times  \cdots \times [-d_s/2, d_s/2] $, and let
 $ \widehat{\chi}_{\fO_{ \bd}}(\bgamma)$ be the Fourier transform of the indicator function $\chi_{ \fO_{\bd}}(\bgamma)$.
It is easy to prove that $ \widehat{\chi}_{\fO_{\bd}}(\b0) =d_1 d_2\cdots d_s$ and
\begin{equation}\label{3.10}
 \widehat{\chi}_{\fO_{\bd}}(\bgamma) = \prod_{i=1}^s \frac{e( d_i \gamma_i/2) -e( -d_i \gamma_i/2)}{2\pi \sqrt{-1}\gamma_i}
    = \prod_{i=1}^s \frac{\sin(\pi d_i  \gamma_i)}{ \pi \gamma_i}  , \quad {\rm for } \quad {\rm Nm}(\bgamma)  \neq 0.
\end{equation}

We fix a nonnegative function $ \omega(\bx),\; \bx\in \RR^s,$ of the
class $ C^\infty $, with a support  inside the unit ball $ |\bx|\leq 1 $, such that
\begin{equation}\label{3}
     \int_{\RR^s }\omega(\bx)d\bx=1.
\end{equation}
We set $\omega_{\tau}(\bx) = \tau^{-s} \omega(\tau^{-1} \bx), \; \tau>0  $, and
\begin{equation}\label{3.0}
 \hat{\omega} (\by)   =        \int_{\RR^s }e(\langle \by,\bx\rangle) \omega(\bx)d\bx.
\end{equation}

 Notice that the Fourier transform $\hat{\omega_{\tau}}(\by)=\hat{\omega}(\tau \by)$ of the
function $ \omega_{\tau}(\by)$ satisfies the bound
\begin{equation}\label{3.1}
|\hat{\omega_{\tau}} (\by)|<c_2(1+\tau |\by|)^{-2s} .
\end{equation}

 {\bf Lemma 4.} {\it There exists a constant $c>0$, such that we have for $  N> c$ 
\begin{equation}\nonumber
   |\cR (\fO_{\btheta \cdot\bN}+\bx,\Gamma ) -  \ddot{\cR} (\fO_{\btheta \cdot\bN}+\bx,\Gamma )| 
    \leq     2^s,
\end{equation} 
where}
\begin{equation}\label{2.34a}
    \ddot{\cR} (\fO_{\btheta \cdot\bN}+\bx,\Gamma) = (\det \Gamma)^{-1}  \sum_{\bgamma\in \Gamma^\bot \setminus  \b0}
   \widehat{\chi}_{\fO_{\btheta \cdot\bN}}(\bgamma)\widehat{\omega}
    (\tau  \bgamma)e(\langle\bgamma,\bx\rangle), \quad \tau =N^{-2}.
\end{equation}
{\bf Proof.} 
Let $\fO^{\pm \tau}_{\btheta \cdot \bN} =  [0, \max(0,\theta_1N_1\pm \tau)) \times \cdots    \times   [0, \max(0,\theta_sN_s\pm \tau))$, and let $\chi_{\fO}(x)$ be the indicator function of $\fO$. 
 We consider the convolutions of the functions $\chi_{\fO^{\pm \tau}_{\bN}}( \bgamma)$  and $\omega_{\tau}(\by)$ : 
\begin{equation}\label{2.14}
    \omega_{\tau} \ast \chi_{\fO^{\pm \tau}_{\bN}}  (\bx) =   \int_{\RR^s}   \omega_{\tau} (\bx-\by)  \chi_{\fO^{\pm \tau}_{\bN}} (\by) d\by . 
\end{equation}
It is obvious that the nonnegative functions (\ref{2.14}) are of class $ C^\infty $ and are compactly supported in 
$\tau$-neighborhoods of the bodies $\fO^{\pm \tau}_{\bN} $, respectively. 
We obtain
\begin{equation}\label{2.16}
     \chi_{\fO_{\btheta \cdot \bN}^{-\tau}} (\bx)   \leq \chi_{\fO_{\btheta \cdot \bN}} (\bx)  \leq  \chi_{\fO_{\btheta \cdot \bN}^{+\tau}} (\bx) , \;\;
     \chi_{\fO_{\btheta \cdot \bN}^{-\tau}} (\bx)  \leq \omega_{\tau}  \ast  \chi_{\fO_{\btheta \cdot \bN} }(\bx)     \leq   \chi_{\fO_{\btheta \cdot \bN}^{+\tau}} (\bx) .
\end{equation}
Replacing $\bx$ by $\bgamma -\bx$ in (\ref{2.16}) and summing these inequalities over  $\bgamma \in\Gamma = \Gamma_M$, we find from (\ref{1.0}), that
\begin{equation}\nonumber
  \cN (\fO^{-\tau}_{\btheta \cdot\bN}+\bx,\Gamma ) \leq    \cN (\fO_{\btheta \cdot \bN}+\bx,\Gamma ) \leq   \cN (\fO^{+\tau}_{\btheta \cdot\bN}+\bx,\Gamma ),
\end{equation}
and
\begin{equation}\nonumber
         \cN (\fO^{-\tau}_{\btheta \cdot\bN}+\bx,\Gamma )    \leq  \dot{\cN} (\fO_{\btheta \cdot\bN}+\bx,\Gamma )   \leq
                        \cN (\fO^{+\tau}_{\btheta \cdot\bN}+\bx,\Gamma ),     
\end{equation}
where
\begin{equation}\label{2.20}
         \dot{\cN} (\fO_{\btheta \cdot\bN}+\bx,\Gamma ) = \sum_{\bgamma\in \Gamma}  \omega_{\tau}  \ast \chi_{\fO_{\btheta \cdot \bN}} (\bgamma - \bx).
\end{equation}
Hence
\begin{equation}\nonumber
   -  \cN (\fO^{+ \tau}_{\btheta \cdot\bN}+\bx,\Gamma ) + \cN (\fO^{- \tau}_{\btheta \cdot\bN}+\bx,\Gamma ) 
\end{equation}
\begin{equation}\nonumber
     \leq  \dot{\cN} (\fO_{\btheta \cdot\bN}+\bx,\Gamma ) - \cN (\fO_{\btheta \cdot\bN}+\bx,\Gamma )  \leq  \cN (\fO^{+ \tau}_{\btheta \cdot\bN}+\bx,\Gamma ) - \cN (\fO^{- \tau}_{\btheta \cdot\bN}+\bx,\Gamma ) .
\end{equation}
Thus
\begin{equation}\label{2.22}
       |\cN (\fO_{\btheta \cdot\bN}+\bx,\Gamma ) -  \dot{\cN} (\fO_{\btheta \cdot\bN}+\bx,\Gamma )| \leq  \cN (\fO^{+ \tau}_{\btheta \cdot\bN}+\bx,\Gamma ) - \cN (\fO^{- \tau}_{\btheta \cdot\bN}+\bx,\Gamma ) .
\end{equation}
Consider the right side of this inequality.
We have that $\fO^{+ \tau}_{\btheta \cdot\bN}  \setminus \fO^{- \tau}_{\btheta \cdot\bN}$ is the union of boxes $\fO^{(i)}, \; i=1,...,2^s-1$, where
\begin{equation}  \nonumber
 \vol(\fO^{(i)}) \leq  \vol( \fO_{\bN}^{+\tau})  -  \vol(\fO_{\bN}^{-\tau})   \leq \prod_{i=1}^s (N_i + \tau) - \prod_{i=1}^s (N_i - \tau) 
\end{equation}
\begin{equation}  \nonumber
 \leq N \Big( \prod_{i=1}^s (1 + \tau) - \prod_{i=1}^s (1 - \tau) \Big) <\ddot{c}_s N \tau =\ddot{c}_s  / N,  , \quad \tau =N^{-2},
\end{equation} 
with some $\ddot{c}_s>0$.  From (\ref{1.01}), we get  $|{\rm Nm}( \bgamma)| \geq C_{M}^s$ for $\bgamma \in \Gamma_{M} \setminus \b0$.
We see $ |{\rm Nm}( \bgamma_1 \\- \bgamma_2  )| \leq  \vol(\fO^{(i)} +\bx) < C_{M}^s $  for $ \bgamma_1, \bgamma_2  \in \fO^{(i)} +\bx $
 and  $  N>\ddot{c}_s / C_{M}^s$.  Therefore, the 
  box $\fO^{(i)}+\bx $ 
contains at most one point of $\Gamma_M$ for $  N>\ddot{c} /C_{M}^s $. By (\ref{2.22}), we obtain
\begin{equation}\label{2.28}
    |\dot{\cN} (\fO_{\btheta \cdot\bN}+\bx,\Gamma ) - \cN (\fO_{\btheta \cdot\bN}+\bx,\Gamma ) |  \leq  2^s -1,  \quad {\rm for}  \quad  N>\ddot{c} / C_{M}^s.
\end{equation}
Let
\begin{equation}\label{2.30}
   \dot{\cR} (\fO_{\btheta \cdot\bN}+\bx,\Gamma ) = \dot{\cN} (\fO_{\btheta \cdot\bN}+\bx,\Gamma ) - \frac{ \vol(\fO_{\btheta \cdot\bN})}{\det \Gamma}   .
\end{equation}  

By  (\ref{2.20}), we have that $ \dot{\cN} (\fO_{\btheta \cdot\bN}+\bx,\Gamma )$ is a periodic
function of $ \bx \in \RR^n $ with the period lattice $
\Gamma$. 
Applying the Poisson summation formula to the series (\ref{2.20}), and bearing in mind that $ \widehat{\omega}_{\tau} (\by) = \widehat{\omega} (\tau \by)$,  
 we obtain from (\ref{2.34a})
\begin{equation}  \nonumber
\dot{\cR} (\fO_{\btheta \cdot\bN}+\bx,\Gamma) = \ddot{\cR} (\fO_{\btheta \cdot\bN}+\bx,\Gamma) .
\end{equation}
Note that (\ref{3.1})
ensures the absolute convergence of the series (\ref{2.34a}) over $\bgamma\in \Gamma^\bot \setminus  \{0\}$.
Using (\ref{1.2}), (\ref{2.28}) and (\ref{2.30}), we get the assertion of Lemma 4.   \qed \\

\setcounter{equation}{0}
\renewcommand{\theequation}{2.4.\arabic{equation}}
{\bf 2.4. Upper bound of the variance  of $\cR (\btheta \cdot \bN  \cdot \KK^s+\bx, \Gamma)  $.} 
Let
\begin{equation}\label{3.23}
       \cA(G) = (\det \Gamma)^{-1}  \sum_{\bgamma \in G}
               \widehat{\chi}_{\fO_{\btheta \cdot\bN}}(\bgamma)
       \hat{\omega}
      (\tau\bgamma)e(\langle\bgamma,\bx\rangle),   
\end{equation}
\begin{equation}\label{3.23a}
      \tilde{ \cA}(G) = (\det \Gamma)^{-1}  \sum_{\bgamma \in G}
       | \widehat{\chi}_{\fO_{\btheta \cdot\bN}}(\bgamma)
              \hat{\omega}
      (\tau\bgamma)|,   
\end{equation}
and let
\begin{equation}\label{3.231}
        \cB(G,\bkappa) =  \sum_{\bgamma \in G} 
       \frac{(\det \Gamma)^{-1} \hat{\omega}(\tau\bgamma)}{(2\pi \sqrt{-1})^s {\rm Nm} (\bgamma)}
      e \Big(\sum_{k=1}^s  \gamma_k(\kappa_k(\theta_k N_k)/2 +x_k) \Big)   .
\end{equation}
We obtain from (\ref{3.10})
\begin{equation}\label{3.232}
        \cA(G)  =\sum_{\kappa_1,...,\kappa_s \in \{-1, 1\}} \kappa_1\kappa_2 \cdots \kappa_s
         \cB(G,\bkappa).   
\end{equation}
Using the Cauchy--Schwartz inequality, we get
\begin{equation}\label{3.233}
      |  \cA(G)|^2  \leq 2^s \sum_{\kappa_1,...,\kappa_s \in \{-1, 1\}} 
      |  \cB(G,\bkappa)|^2.
\end{equation}
By (\ref{3.20}) and (\ref{2.34a}), we see that
\begin{equation}\label{3.24}
   \ddot{\cR} (\fO_{\btheta \cdot\bN}+\bx,\Gamma ) =    \cA(G^{(1)} ) +...+  \cA(G^{(5)})   .
\end{equation}
Let 
\begin{equation}\label{3.24a}
       h(\bgamma) =
       \frac{(\det \Gamma)^{-1} \hat{\omega}(\tau\bgamma)}{(2\pi \sqrt{-1})^s {\rm Nm} (\bgamma)}
      e\Big(\sum_{k=1}^s  \gamma_k  x_k \Big).
\end{equation}
It is easy to see that 
\begin{equation}\label{3.24b}
    \cB(G,\bkappa) = \sum_{\bgamma \in G}    h(\bgamma)   e\Big(\sum_{k=1}^s \kappa_k \gamma_k\theta_k N_k /2\Big),
\end{equation} 
and
\begin{equation}\label{3.24c}
    \cA(G) = \sum_{\kappa_1,...,\kappa_s \in \{-1, 1\}} \kappa_1\kappa_2 \cdots \kappa_s
    \sum_{\bgamma \in G}    h(\bgamma)   e\Big(\sum_{k=1}^s \kappa_k \gamma_k\theta_k N_k /2\Big).
\end{equation} \\

{\bf Lemma 5.} {\it With notations as above }
\begin{equation}\nonumber
       \cA(G^{(2)})  =O(1/N)   .
\end{equation}
\\
{\bf Proof.} By (\ref{3.231}) and (\ref{3.1}) we have that
\begin{equation}\label{3.25}
       | \cB(G^{(2)},\bkappa)|  \leq c_2  \sum_{\bgamma \in G^{(2)}} 
       \frac{(\det \Gamma)^{-1} (\tau |\bgamma|)^{-2s}}{(2\pi) ^s |{\rm Nm} (\bgamma)|}, \quad \tau=N^{-2}.
\end{equation}
Notice that for every lattice $ \fL\in \RR^s $, one has
the bound (see, e.g., [GL] p. 141, 142)
\begin{equation}\nonumber
    \#\{\bgamma \in \fL \; : \; j \leq  | \bgamma | \leq j+1\}  =O(j^{s-1})   .
\end{equation}
Hence
\begin{equation}\nonumber
  \sum_{\bgamma\in \Gamma^\bot: |\bgamma| \geq N^5 }  |\bgamma|^{-2s}  \leq  \sum_{j \geq N^5}  \sum_{\bgamma\in \Gamma^\bot: |\bgamma|\in [j,j+1) }|\bgamma|^{-2s}  = O\Big(\sum_{j \geq N^5} j^{-s-1}\Big) =O(N^{-5s}) .
\end{equation}  
By  (\ref{1.01}), (\ref{3.21}) and (\ref{3.25}), we obtain
\begin{equation}\nonumber
       | \cB(G^{(2)},\bkappa)|  \leq c_2  (C_M^s\det \Gamma(2\pi) ^s )^{-1}  \sum_{\bgamma \in \Gamma^\bot: 
       |\bgamma| \geq N^5}         N^{4s}|\bgamma|^{-2s}   =O(N^{-s}).
\end{equation}
Using (\ref{3.232}),  we get the assertion of Lemma 5. \qed

We consider the probability space $([0,1]^s, \lambda , B([0,1]^s))$ with Lebesgue's measure $\lambda$. 
Hence, we have the following formula for the expectation:
\begin{equation}\label{3.26}
      {\bf E} [f(\btheta)]  = \int_{[0,1]^s} f(\btheta) d \btheta .              
\end{equation}

{\bf Lemma 6.} {\it Let  $\bgamma^{(i)} \in \Gamma^{\bot}, \; i=1,2, \;   \bgamma^{(1)}  \neq \bgamma^{(2)}$. Then   }
\begin{equation}\nonumber
      |{\bf E} [ e(\langle\bgamma^{(1)}-\bgamma^{(2)},\btheta \cdot \bN \rangle/2 +\beta)] |  \leq    \frac{1}{\pi^s C_M^s  N } .
\end{equation}

{\bf Proof.} Bearing in mind that 
\begin{equation}\label{3.262}
   \Big|\int_0^1 e(x  z )dx \Big| =\Big|\frac{e(z)-1}{2\pi z}\Big| \leq \frac{1}{|\pi z|},  \quad {\rm with} \quad z \neq 0,
\end{equation}      
   (\ref{1.01})     and that $N_1 \cdots N_s=N$,           we have 
\begin{equation}\nonumber
      |{\bf E} [ e(\langle\bgamma^{(1)}-\bgamma^{(2)},\btheta \cdot \bN \rangle/2 +\beta)]| 
      \leq \frac{1}{\pi^s N |{\rm Nm} (\bgamma^{(1)}-\bgamma^{(2)})|}     \leq    \frac{1}{\pi^s C_M^s N }. \qquad \qed
\end{equation} \\

{\bf Lemma 7.}  {\it With notations as above }
\begin{equation}\label{3.27}
      {\bf E} [ |  \cA(G^{(1)})|^2]  =O(n^{s-1})  , \qquad  {\bf E} [ |  \cA(G^{(3)})|^2]  =O(n^{s-3/2}) ,
\end{equation}
and
\begin{equation}\label{3.27b}
      |  \cA(\dot{G})| \leq    \tilde{\cA}(\dot{G})=O(n^{s})    ,  \qquad {\rm with}   \quad G \subseteq G^{(1)}  .
\end{equation}
\\

{\bf Proof.} By (\ref{3.231}), (\ref{3.233}) and the Cauchy--Schwartz inequality, we obtain that
\begin{equation} \nonumber
      {\bf E} [ | \cA(G)|^2]  \leq  2^s \sum_{\kappa_1,...,\kappa_s \in \{-1, 1\}} \Big( \sum_{\bgamma \in G}
       \frac{(\det \Gamma)^{-2} |\hat{\omega}(\tau\bgamma)|^2}{(2\pi)^{2s} |{\rm Nm} (\bgamma)|^2} 
\end{equation}
\begin{equation} \nonumber 
      + \sum_{\bgamma^{(1)},\bgamma^{(2)}  \in G, \;\; \bgamma^{(1)} \neq \bgamma^{(2)}} 
       \frac{(\det \Gamma)^{-2} |\hat{\omega}(\tau\bgamma^{(1)}) ||\hat{\omega}(\tau\bgamma^{(2)}) |}
       {(2\pi)^{2s} |{\rm Nm} (\bgamma^{(1)}) ||{\rm Nm} (\bgamma^{(2)}) |} 
      |{\bf E} [ e(\langle\bgamma^{(1)}-\bgamma^{(2)},\bkappa \cdot\btheta \cdot N \rangle/2 +\beta)]  | \Big)
\end{equation}
with $\beta = \sum_{k=1}^s (\gamma_k^{(1)} -\gamma_k^{(2)} )x_k $.\\

Applying Lemma 6, we get
\begin{equation}\label{3.271}
      {\bf E} [ |  \cA(G)|^2]  \leq  2^{2s}\sum_{\mu \in [1,s]} S_{1, \mu}(G)+  \frac{2^{2s}}{\pi^s C_M^s  N } S_2(G,G),
\end{equation}
with
\begin{equation}\label{3.300}
   S_{1, \mu}(G) =   \sum_{\bgamma \in G, |\gamma_i| \leq |\gamma_{\mu}|, i=1,....,s} 
       \frac{(\det \Gamma)^{-2} |\hat{\omega}(\tau\bgamma)|^2}{(2\pi)^{2s} |{\rm Nm} (\bgamma)|^2} 
\end{equation}
and
\begin{equation} \label{3.305}
    S_2(\dot{G}, \ddot{G})=   \sum_{\bgamma^{(1)} \in \dot{G}, \bgamma^{(2)}  \in \ddot{G}, \;\; \bgamma^{(1)} \neq \bgamma^{(2)}} 
             \frac{(\det \Gamma)^{-2} |\hat{\omega}(\tau\bgamma^{(1)})||\hat{\omega}(\tau\bgamma^{(2)})|}
       {(2\pi)^{2s} |{\rm Nm} (\bgamma^{(1)}) ||{\rm Nm} (\bgamma^{(2)}) |}. 
\end{equation}
We fix $\mu \in [1,s]$, and we consider $ S_{1, \mu}(G)$. 
 Let $\bgamma \in     B(m, \mu,\bkappa,j)$. 
According to (\ref{2.102}) and (\ref{2.121}),   we have
\begin{equation} \label{3.31}
   |{\rm Nm}( \bgamma)|= (j+z_1)C_M^s 2^{z_2(s-1)},   \quad  {\rm and} \quad \#\Gamma^{\bot}   \cap  B(\bm, \mu,\bkappa,j) \leq 1
\end{equation}  
 with $z_1,z_2 \in [0,1]$. 
By (\ref{3.21}) and (\ref{3.31}), we obtain
\begin{equation}\nonumber
   j+1  \geq    n^{1/2}(C_M^s 2^{s-1})^{-1}   \quad {\rm for} \quad  \bgamma  \in G^{(3)}.
\end{equation}
Hence 
\begin{equation}\label{3.33}
       \sum_{j \geq 0}  \sum_{\bgamma \in G^{(3)} \cap  B(m, \mu,\bkappa,j) }
        \frac{1}{|{\rm Nm} (\bgamma)|^2}  =O\Big(\sum_{j \geq n^{1/2}}   1/j^2 \Big) =O(n^{-1/2}),
\end{equation}
\begin{equation}\label{3.34}
       \sum_{j \geq 0}  \sum_{\bgamma \in G \cap  B(m,\mu, \bkappa,j) }
        \frac{1}{|{\rm Nm} (\bgamma)|^2}   =O\Big(\sum_{j \geq 1}   1/j^2\Big) =O(1), \quad  {\rm for} \quad G \subseteq G^{(1)}.
\end{equation}
Bearing in mind that $ {\rm Nm} (\bgamma) \leq N^{5s}$ for $\bgamma \in G^{(1)} \cup G^{(3)}$ and $n=[\log_2 N]+1$, we get from 
(\ref{3.21})  and (\ref{3.31}) that
\begin{equation}\label{3.35}
       \sum_{j \geq 0}  \sum_{\bgamma \in (G^{(1)} \cup G^{(3)}) \cap  B(m,\mu, \bkappa,j) }
        \frac{1}{|{\rm Nm} (\bgamma)|}  = O\Big(\sum_{1 \leq j \leq N^{5s}}   1/j \Big)  =O(n).
\end{equation}
By (\ref{2.102}), (\ref{3.21}) and (\ref{1.01}),  we have for $\bgamma \in (G^{(1)} \cup G^{(3)})\cap  B(m,\mu, \bkappa,j) $ that
\begin{equation}\nonumber
 \log_2 |\gamma_i| \in [m_i, m_i +1), \quad   i \in [1,s], \; i \neq \nu(\mu),  \quad |\gamma_i| \leq  N^5, \quad  |{\rm Nm} (\bgamma)| \geq C_M^s ,
\end{equation}
\begin{equation}\nonumber
   C_M^s N^{-5(s-1)}  \leq  |\gamma_i|,   \quad {\rm and} \quad   s\log_2C_M  -5(s-1)\log_2 N  \leq   \log_2 |\gamma_i|  \leq   5\log_2 N, \; 1 \leq i \leq s . 
\end{equation}
Therefore
\begin{equation}\label{3.36}
   \quad m_i \in  J, \quad i \in [1,s], \; i \neq \nu(\mu) \quad {\rm with} \quad J= [s\log_2C_M  -5(s-1)n  ,5n] .
\end{equation}
From (\ref{2.102a}), (\ref{3.1})  and (\ref{3.300}), we derive
\begin{equation} \label{3.36a}
   S_{1, \mu} (G) \leq  \sum_{\kappa_1,...,\kappa_s \in \{-1, 1\}}  
          \sum_{\bm \in  \ZZ^s  \atop m_1+...m_s=0}  \sum_{j  \geq 0}
       \sum_{\bgamma \in G \cap  B(\bm,\mu, \bkappa,j) \atop   |\gamma_i| \leq |\gamma_{\mu}|, i=1,...,s  }
        \frac{ (\det \Gamma)^{-2} c_2^2}{(2\pi)^{2s} |{\rm Nm} (\bgamma)|^2} .
\end{equation}
Hence, we obtain for $i=1,3$ that
\begin{equation} \nonumber
   S_{1, \mu} (G^{(i)}) \leq  \sum_{\kappa_1,...,\kappa_s \in \{-1, 1\}}  
          \sum_{m_k  \in J,  k \in [1,s], \; k \neq \nu(\mu) \atop m_1+...+m_s=0}  \sum_{j  \geq 0}
       \sum_{\bgamma \in G^{(i)} \cap  B(\bm,\mu, \bkappa,j) }
        \frac{ (\det \Gamma)^{-2} c_2^2}{(2\pi)^{2s} |{\rm Nm} (\bgamma)|^2} .
\end{equation}
Applying (\ref{3.33}), (\ref{3.34}) and (\ref{3.36}), we get that
\begin{equation}\label{3.37}
   S_{1,\mu}(G^{(1)}) = O(n^{s-1}),  \quad \quad {\rm and}  \quad \quad      S_{1,\mu}(G^{(3)}) = O(n^{s-1-1/2}).
\end{equation}
Analogously, we have from (\ref{3.305}) and (\ref{3.35}) that   for $\dot{G}, \ddot{G} \subset G^{(1)} \cup G^{(3)}$
\begin{equation}\label{3.38}
    S_2(\dot{G}, \ddot{G}) =O\Big(\Big( \sum_{\kappa_1,...,\kappa_s \in \{-1, 1\}}  
          \sum_{m_k  \in J,  k \in [1,s] \atop  k \neq \nu(\mu), m_1+...+m_s=0}  \sum_{j  \geq 0}
       \sum_{\bgamma \in G^{(1)} \cup G^{(3)} \cap  B(\bm,\mu, \bkappa,j) }
        \frac{ 1}{{\rm Nm} (\bgamma)}   \Big)^2\Big)    
    =O(n^{2s}). 
\end{equation}
According to (\ref{3.271}), we obtain  (\ref{3.27}). 
By (\ref{3.23a}) and (\ref{3.10}), we have that
$$
   \tilde{ \cA}(G) \leq  (\det \Gamma)^{-1} c_2 \sum_{\bgamma \in G}
       1/|{\rm Nm} (\bgamma)|. 
$$
Now using (\ref{3.35}),  similarly to (\ref{3.36a})-(\ref{3.37}), we get  (\ref{3.27b}). 
Hence, Lemma 7 is proved. \qed \\

{\bf Lemma 8.} {\it With notations as above 
\begin{equation}\nonumber
      S_{1, \mu}( \tilde{G}) =O( n^{s-11/9} ) ,
       \quad {\it with}  \quad   \tilde{G}  = \dot{G}_0 \bigcup_{1 \leq i \leq n^{5/9 } +1} \ddot{G}_i , 
       \quad \mu=1,...,s.  
\end{equation} }

{\bf Proof.} Let 
$
   \bgamma  \in   \ddot{G}_{i}\cap  B(\bm, \mu,\bkappa,j).
$
By (\ref{2.102}), we have that $\log_2 |\gamma_k| \in [m_k, m_k +1) $     
   with  $ k \in [1,s], \; k \neq \nu(\mu) $.
From (\ref{3.56}) and (\ref{1.01}),  we derive for $  |\gamma_i| \leq |\gamma_{\mu}|, i=1,....,s$ that
\begin{equation}\nonumber
   \log_2 |\gamma_k|  \leq   (i+1) n^{4/9},    \quad 
  \log_2 |\gamma_k|    \geq  s\log_2C_M   - \sum_{j \in [1,s], \; j \neq k}\log_2 |\gamma_j|, \;k =1,...,s  ,
\end{equation}
and 
\begin{equation}\nonumber
  \log_2 |\gamma_{\mu}|  > (i+1) n^{4/9} - n^{2/9}.
\end{equation}
Therefore
\begin{equation}\nonumber
        m_{\mu} \in J_1,      \quad {\rm with}  \quad     J_1 =( (i+1) n^{4/9} - n^{2/9}-1,(i+1) n^{4/9}   ] , \quad \#J_1  \leq n^{2/9}+2,
\end{equation}
and
\begin{equation}\nonumber
        m_{k} \in J_2      \quad {\rm with}  \quad     J_2 =( -(s-1)(i+1) n^{4/9} +s\log_2 C_M -1  ,(i+1) n^{4/9}   ] ,
\end{equation}
  $    k \in [1,s], \; k \neq \nu(\mu), \mu , \;\#J_2 =O((i+1) n^{4/9}),\; \; i=1,2,...$.
By (\ref{3.36a}), we get that
\begin{equation} \nonumber
   S_1(\ddot{G}_{i}) \leq   \sum_{\kappa_1,...,\kappa_s \in \{-1, 1\}}  
    \sum_{m_{\mu} \in J_1} \;\;
          \sum_{m_k \in J_2\atop  k\neq \mu, \nu(\mu)}  \; \sum_{j  \geq 0} \;
       \sum_{\bgamma \in \ddot{G}_{i} \cap  B(\bm,\mu, \bkappa,j) \atop   |\gamma_i| \leq |\gamma_{\mu}|, i=1,....,s}
        \frac{(\det \Gamma)^{-2} c_2^2  }{(2\pi)^{2s} |{\rm Nm} (\bgamma)|^2}  .
\end{equation}
Using (\ref{3.34}), we obtain
\begin{equation}\nonumber
   S_{1,\mu}(\ddot{G}_i) =O(\#J_1\#J_2^{s-2})=O( i^{s-2} n^{((s-2)4+2)/9} ).
\end{equation}
Similarly we get that $ S_{1,\mu}(\dot{G}_0)  = O(  n^{4(s-1)/9} ) $.
Now from (\ref{3.300})  and (\ref{3.56}), we have
\begin{equation}\nonumber
   S_{1,\mu}(\tilde{G}) = S_{1,\mu}(\dot{G}_0) +   \sum_{1 \leq i \leq n^{5/9 } +1}S_{1,\mu}(\ddot{G}_i) = 
    O\Big( \sum_{1 \leq i \leq n^{5/9 }}  i^{s-2} n^{((s-2)4+2)/9} \Big) =O( n^{s-1-2/9} ) .
\end{equation}
Hence, Lemma 8 is proved. \qed    
\\

{\bf Lemma 9.} {\it With notations as above }
\begin{equation}\nonumber
      {\bf E} [ |  \cA(G^{(5)})|^2]  =O(n^{s-2}  (\ln n)^4)   .
\end{equation}
{\bf Proof.} Let
\begin{equation}\label{3.39a}
       G^{(5,\mu)} = \{\bgamma \in   G^{(5)} \; : \; |N_\mu \gamma_\mu| \leq 2^{(\ln n)^4}  
               \quad {\rm and} \quad            |N_j \gamma_j| > 2^{(\ln n)^4}  \quad {\rm for}   \quad  j<\mu      \} .                
\end{equation}
By (\ref{3.22}), we have that
\begin{equation}\nonumber
        G^{(5)} = \bigcup_{\mu \in [1,s]} G^{(5,\mu)}, \qquad  G^{(5,\mu)} \cap G^{(5,j)} =\emptyset \;\;   {\rm for} \;\;  \mu \neq j.
\end{equation}
Similarly to (\ref{3.271})-(\ref{3.305}), using   the Cauchy--Schwartz inequality,  we obtain  from  (\ref{3.10}) and (\ref{3.23}) that
\begin{equation}\nonumber
      {\bf E}[  |  \cA(G^{(5)})|^2]  \leq s2^{s-1} \sum_{\mu \in [1,s]} \sum_{\kappa_j \in \{-1, 1\} \atop j \in [1,s],  j\neq \mu} 
       \sum_{\bgamma^{(1)},\bgamma^{(2)}  \in G^{(5,\mu)}} \dot{\psi}(\bgamma^{(1)},\bgamma^{(2)}),     
\end{equation}
with
\begin{equation}\label{3.40ab}
      \dot{\psi}(\bgamma^{(1)},\bgamma^{(2)}) = 
       \frac{(\det \Gamma)^{-2} |\hat{\omega}(\tau\bgamma^{(1)}) ||\hat{\omega}(\tau\bgamma^{(2)}) |}
       {(2\pi)^{2s} |{\rm Nm} (\bgamma^{(1)}) ||{\rm Nm} (\bgamma^{(2)}) |}   \tilde{\psi}(\bgamma^{(1)},\bgamma^{(2)})
\end{equation}
and
\begin{equation}\label{3.40b}
     \tilde{\psi}(\bgamma^{(1)},\bgamma^{(2)}) =
       \Big|{\bf E} \Big[ \sin(\pi\theta_\mu N_\mu\gamma_\mu^{(1)}) 
\end{equation}
\begin{equation}\nonumber
     \times  \sin(\pi \theta_\mu N_\mu  \gamma_\mu^{(2)})  e\Big(  \sum_{1 \leq j \leq s, \; j \neq \mu}  (\gamma_j^{(1)} -\gamma_j^{(2)} )
     \theta_j N_j  \kappa_j/2   \Big)   \Big]  \Big|.
\end{equation}
Hence
\begin{equation}\label{3.40}
      {\bf E}[  |  \cA(G^{(5)})|^2]  \leq s2^s \sum_{\mu \in [1,s]} \sum_{\kappa_j \in \{-1, 1\}, j \in [1,s], \; j\neq \mu} 
   (\dot{S}_1(\mu)+\dot{S}_2(\mu))
\end{equation}
with
\begin{equation}\label{3.41}
    \dot{S}_1(\mu) =\sum_{\bgamma \in G^{(5,\mu)} }   \dot{\psi}(\bgamma,\bgamma)  \quad {\rm and} \quad
     \dot{S}_2(\mu) = \sum_{\bgamma^{(1)},\bgamma^{(2)}  \in G^{(5,\mu)},  \atop \bgamma^{(1)} \neq \bgamma^{(2)}}  \dot{\psi}(\bgamma^{(1)},\bgamma^{(2)}). 
\end{equation}
Bearing in mind that $ |\sin(x)| \leq |x|$, we derive from   (\ref{3.40b}) that
\begin{equation}\label{3.41c}
        \tilde{\psi}(\bgamma,\bgamma)    \leq  \min(1,|2\pi N_\mu \gamma_\mu|^2).
\end{equation}
 Consider $  \dot{S}_2(\mu)$.   By   (\ref{3.262}),    we get for $ \bgamma^{(1)} \neq \bgamma^{(2)} $   that
\begin{align} \nonumber
        \tilde{\psi}(\bgamma^{(1)},\bgamma^{(2)})  & \leq  \prod_{j \in [1,s], j \neq \mu} \frac{1}{\pi N_j |\gamma_j^{(1)}-\gamma_j^{(2)}|} \\
    & =  \frac{N_\mu|\gamma_{\mu}^{(1)}-\gamma_{\mu}^{(2)}|}{\pi^{s-1} N |{\rm Nm} (\bgamma^{(1)}-\bgamma^{(2)})|}  \leq 
    \frac{2^{(\ln n)^4+1}}{\pi^{s-1} C_M^s  N }. \nonumber 
\end{align}
According to  (\ref{2.90}), (\ref{3.22}) and (\ref{1.01}), we have   $G^{(5)}  \subset  \GG(s\log_2 C_M-(s-1)(n+1), n+1)$.
 Using Lemma 2, we obtain
 $\#G^{(5)}  =O(n^{s})$. Applying  (\ref{3.40ab}) and (\ref{3.41}), we get
\begin{equation}\label{3.41b}
      \dot{S}_2(\mu) =O(   n^{2s}  N^{-1}2^{(\ln n)^4+1}) =O(1).
\end{equation}
Now we fix $\mu \in [1,s]$, and we consider $  \dot{S}_1(\mu)$.  Let
$$
   \bgamma \in \Gamma^{\bot}\cap  B(\bm, \mu,\bkappa^{(1)},j).
$$
According to  (\ref{2.102}) and (\ref{3.39a}),   we have that
\begin{equation}\nonumber
     \log_2 |N_{\mu} \gamma_{\mu}| =  \log_2 N_{\mu} +m_{\mu} +z_1 \leq (\ln n)^4,  \quad z_2 \in [0,1) .
\end{equation}
Hence
\begin{equation}\nonumber
     m_{\mu} \in \dot{J},   \quad {\rm with} \quad   \dot{J} =(-\infty,  (\ln n)^4   -  \log_2 N_{\mu}  ].
\end{equation}
By (\ref{3.31}) and (\ref{1.01}), we obtain that   $ C_M^s \max(1,j) \leq  {\rm Nm} (\bgamma)   \leq  (j+1) C_M^s 2^{s-1}$ and
\begin{equation} \nonumber
    \sum_{m_{\mu} \in \dot{J}}       \sum_{j \geq 0}  \sum_{\bgamma \in G^{(5,\mu)} \cap  B(\bm, \mu,\bkappa^{(1)},j) }
        \frac{\min(1,|N_\mu\gamma_\mu|^2)}{|{\rm Nm} (\bgamma)|^2} 
        =O\Big( \sum_{m_{\mu} \in \dot{J}}       \sum_{j \geq 1} \frac{\min(1,N_{\mu}^2 2^{2m_\mu})}{j^2}\Big) 
\end{equation}
\begin{equation} \label{3.43}
        =O\Big(  \sum_{m_{\mu} \leq  -  \log_2 N_{\mu}}       \sum_{j \geq 1} \frac{N_{\mu}^2 2^{2m_\mu}}{j^2}   +
         \sum_{m_{\mu} \in  [ -  \log_2 N_{\mu} , (\ln n)^4   -  \log_2 N_{\mu}]}       \sum_{j \geq 1} \frac{1}{j^2}\Big) =O((\ln n)^4).
\end{equation}
Bearing in mind (\ref{2.102}) and that $|\bgamma| \leq N$, we have for $k \neq \mu, \nu(\mu)$ that
\begin{equation} \nonumber
  m_k +z_1=   \log_2 | \gamma_{k}| =  \log_2 {\rm Nm}(\bgamma)  - \sum_{j \in [1,s], j \neq k}  \log_2 | \gamma_{j}| 
       \geq \log_2 C_M^s  -(s-1)(n+1), 
\end{equation}
with $z_1 \in [0,1) $  and  
\begin{equation}\label{3.43a}
     m_k \in  \ddot{J}   \quad {\rm with} \quad    \ddot{J} =[\log_2 C_M^s -(s-1)(n+1) -1,  n+1].
\end{equation} 
By (\ref{2.102a}), (\ref{3.40ab}), (\ref{3.41}),  (\ref{3.41c}) and  (\ref{3.1}), we get
\begin{equation} \nonumber
   \dot{S}_1(\mu) \leq     \bigcup_{\bkappa^{(1)} \in \{-1, 1\}^s}   \bigcup_{\bm \in \ZZ^s, \atop m_1+...+m_s=0} 
\bigcup_{j \geq 0}   \sum_{\bgamma \in G^{(5,\mu)} \cap  B(\bm,\mu, \bkappa^{(1)},j) }   \dot{\psi}(\bgamma,\bgamma)
\end{equation}
\begin{equation} \nonumber
   \leq    \bigcup_{\bkappa^{(1)} \in \{-1, 1\}^s}  \sum_{m_{\mu} \in \dot{J}}
          \sum_{m_k \in \ddot{J}, \atop k \in [1,s], k \neq \mu,\nu(\mu) } \;  \sum_{j  \geq 0} \;
       \sum_{\bgamma \in G^{(5,\mu)} \cap  B(\bm,\mu, \bkappa^{(1)},j) }
        \frac{(\det \Gamma)^{-2} c_2^2 \min(1, |2\pi N_\mu\gamma_\mu|^2)}{(2\pi)^{2s} {|\rm Nm} (\bgamma)|^2}  .
\end{equation}
Applying  (\ref{3.43}) and (\ref{3.43a}), we derive
\begin{equation}\nonumber
   \dot{S}_1(\mu) = O(\#\ddot{J}^{s-2}  (\ln n)^4) = O(n^{s-2}  (\ln n)^4). 
\end{equation}
By (\ref{3.40}) and (\ref{3.41b}), Lemma 9 is proved. \qed 
\\

 {\bf Lemma 10.}   {\it With notations as above  }
\begin{equation}\label{3.45}
       \cA(G^{(4)}) =O(1)   .
\end{equation}
{\bf Proof.} By (\ref{3.1}), (\ref{3.10}) and (\ref{3.23}), we have
\begin{equation}\nonumber
     | \cA(G^{(4)})| \leq  (\det \Gamma)^{-1} c_2 \sum_{\bgamma \in G^{(4)}}
        \prod_{i=1}^s \frac{|\sin (\pi(\theta_i N_i  \gamma_i)|}{2\pi |\gamma_i| } .
\end{equation}
From (\ref{3.21}), we get for  $\bgamma \in  G^{(4)}$  that $ |\bgamma|>N$,
\begin{equation}\nonumber
  \exists \nu \in [1,s] \quad 
   \quad {\rm with } \quad \quad \log_2 (|\gamma_\nu|) \geq  \log_2 (|\bgamma|/s) \geq n - 1 -\log_2 s   ,
\end{equation}
and
\begin{equation}\nonumber
     \log_2 (|\gamma_1|) +...+ \log_2 (|\gamma_s|) \leq 1/2 \log_2 n, \qquad n=[\log_2 N] +1  .
\end{equation}
Hence,  there exists   $\mu \in \{1,...,s  \} \setminus \{ \nu  \}$ with
\begin{equation}\nonumber
  \log_2 (|\gamma_\mu|) \leq   \frac{1}{s-1} \sum_{i \in [1,s], i \neq \nu} \log_2 (|\gamma_i|)  =
                    \frac{1}{s-1} \Big(\sum_{i \in [1,s]} \log_2 (|\gamma_i|) -  \log_2 (|\gamma_{\nu}|) \Big),
\end{equation}
\begin{equation}\nonumber
    \leq ( - n +1 +\log_2 s  + 1/2 \log_2 n )/(s-1)  , \qquad {\rm and}  \quad |\gamma_\mu| \leq 4 N^{ -\frac{1}{s-1} }n^{1/2}.
\end{equation}
 Bearing in mind that $N_\mu N^{-1/s} \in [1/c_0,c_0]$, we obtain
\begin{equation}\nonumber
   |\sin (\pi(\theta_\mu N_\mu \gamma_\mu))| \leq \pi|\theta_\mu N_\mu \gamma_\mu | 
   \leq \pi |N_\mu \gamma_\mu|  =O(N^{\frac{1}{s} -\frac{1}{s-1} }n^{1/2}) =O(n^{-s}).
\end{equation}
Therefore
\begin{equation} \nonumber
     |  \cA(G^{(4)})|=O\Big( n^{-s} \sum_{\bgamma \in G^{(4)}} 1/ | {\rm Nm}(\bgamma)| \Big) = O( n^{-s} \# G^{(4)}).
\end{equation}
Taking into account that $G^{(4)} \in \GG(0,10n) $,  we get from Lemma 2   the assertion of  Lemma~10. 
\hfill  \qed  \\

{\bf Lemma 11.} { \it There exists a real $w_2 >0$ such that
\begin{equation}\label{3.54}
     {\bf E} [ (\cR (\btheta \cdot \bN  \cdot \KK^s+\bx, \Gamma) )^2 ]  \leq w_2 n^{s-1} ,
\end{equation}
\begin{equation}\label{3.55}
     {\bf E} [| \cR (\btheta \cdot \bN  \cdot \KK^s+\bx, \Gamma) -  \cA(G^{(1)}) |^2]  =O(n^{s-1-1/2}),
\end{equation}
and}
\begin{equation}\label{4.50a}
     {\bf E} \Big[ (\cR (\btheta \cdot \bN  \cdot \KK^s+\bx, \Gamma) - \sum_{i\in [1,n^{5/9}]}       \cA(\dot{G}_i) )^2 \Big] =O ( n^{s-1-2/9}) .
\end{equation}

{\bf Proof.} 
By  (\ref{3.24}) and Lemma 4, we get
\begin{equation}\nonumber
 |\cR (\btheta \cdot \bN  \cdot \KK^s+\bx, \Gamma) -   \sum_{i=1}^5 \cA(G^{(i)}) | 
 \leq   2^s.  
\end{equation} 
It is easy to see that
\begin{equation}\nonumber
 |\cR (\btheta \cdot \bN  \cdot \KK^s+\bx, \Gamma) -   \cA(G^{(1)}) | \leq 
 \sum_{i=2}^5  |  \cA(G^{(i)})|    +  2^s.
\end{equation}  
Using the Cauchy--Schwartz inequality, we obatain
\begin{equation}\nonumber
  {\bf E} [ |\cR (\btheta \cdot \bN  \cdot \KK^s+\bx, \Gamma) -   \cA(G^{(1)}) |^2] 
 \leq    5 \Big(  \sum_{i=2}^5  {\bf E} [|  \cA(G^{(i)})|^2]  +  2^{2s} \Big).
\end{equation} 
Applying Lemma 5 - Lemma 10, we have (\ref{3.55}). By Lemma 7, the triangle inequality and the Cauchy--Schwartz inequality, we get (\ref{3.54}).

Now consider the statement (\ref{4.50a}).
From (\ref{3.57}) and (\ref{3.23}) we obtain, that
\begin{equation}\label{4.59}
         \cA(G^{(1)})  = \sum_{i=1}^{[n^{5/9}]}   \cA(\dot{G}_i)  +      \cA(\tilde{G}), \quad {\rm with}   \quad    
          \tilde{G} =    \dot{G}_0 \bigcup_{i=1}^{[n^{5/9}]} \ddot{G}_i.
\end{equation}
According to (\ref{3.271}), we have  
\begin{equation} \nonumber
      {\bf E} [ |  \cA(\tilde{G})|^2]  \leq  2^{2s}\sum_{\mu \in [1,s]} S_{1, \mu}(\tilde{G})+  \frac{2^{2s}}{\pi^s C_M^s  N }
       S_2(\tilde{G},\tilde{G}).
\end{equation}
Using Lemma 8 and (\ref{3.38}), we derive
\begin{equation} \nonumber
     {\bf E} [  |\cA(\tilde{G})|^2]     =O ( n^{s-1-2/9}) .
\end{equation}
From (\ref{3.55}) and the triangle inequality, we get  (\ref{4.50a}). 
 Therefore, Lemma 11 is proved. \\ \hfill  \qed  
\\

{\bf 2.5. Lower bound of variance of $\cR (\btheta \cdot \bN  \cdot \KK^s+\bx, \Gamma)  $.}\\

\setcounter{equation}{0}
\renewcommand{\theequation}{2.5.\arabic{equation}}
 {\bf Lemma 12.} {\it There exist reals $w_3, c_5>0$ such that  for $N>c_5$ and for all $\bx$}
\begin{equation}\label{4.0}
  \check{ \psi} :=   \sum_{\bgamma \in G^{(1)}} 
       \frac{(\det \Gamma)^{-2} }{(2\pi)^{2s} |{\rm Nm} (\bgamma)|^2}  
        {\bf E} \Big[    \prod_{i=1}^s  \sin^2 (\pi\theta_i N_i\gamma_i)\Big]  2 \cos^2(2\pi\langle\bgamma,\bx\rangle  ) \geq w_3 n^{s-1} .
\end{equation} \\

{\bf Proof.}   Let $\bm \in \ZZ^s, m_1+ \cdots +m_s=0$, $q =2+[\det \Gamma^\bot/C_M^s]$, and
\begin{equation}\nonumber
   \cD(\bm) =\prod_{i=1}^s [-q^{m_i},q^{m_i} ] \times [-\det \Gamma^\bot q^{m_{s}},\det \Gamma^\bot q^{m_{s}} ].
\end{equation}
According to Minkowsky's theorem,  there exists $\bgamma(\bm) \in \Gamma^{\bot} \setminus \{0\}$ with  
$\bgamma(\bm) \in  \cD(\bm)$. We see that
\begin{equation}\label{4.0a}
   | {\rm Nm} (\bgamma(\bm))  | \leq  \det \Gamma^\bot.  
\end{equation}
Suppose $|\gamma(\bm)_i| \leq q^{m_i-1} $ for some $i \in [1,s-1]$.
By (\ref{1.01}) we get
\begin{equation}\nonumber
   C_M^s  \leq | {\rm Nm} (\bgamma(\bm))  | \leq  \det \Gamma^\bot/q <    C_M^s .  
\end{equation}
We arrive at a contradiction. Hence
\begin{equation}\label{4.0b}
 |\gamma(\bm)_i \; | \; \in (q^{m_i-1},q^{m_i}]  \; {\rm for }  \; i \in [1,s-1],
    \;  {\rm and }  \; \bgamma(\bm^{(1)}) \neq \bgamma(\bm^{(2)}),  
\end{equation}
 for  $\bm^{(1)} \neq \bm^{(2)}$. 
Let
\begin{equation}\label{4.0c}
   \bar{G} =\{\bgamma(\bm) \; | \;  -n/(4s) \log_q 2  \leq m_i \leq -2s , \; i=1,...,s-1  \}.
\end{equation}
We see for sufficiently large $N$ that
\begin{equation}\label{4.3}
   \#\bar{G} \geq  n^{s-1} ((5s)^{-1} \log_q 2)^{s-1}.
\end{equation}  
By  (\ref{2.81})  $N_iN^{1/s} \in [1/c_0,c_0]$.  From (\ref{4.0c}),  we obtain for sufficiently large $N$ that
\begin{equation}\nonumber
   N_i 2^{-2s} \geq |N_i\gamma_i | \geq c_0^{-1} 2^{n/s -n/(4s) \log_q 2 -2 }  \geq 2^{ln^4n}, \quad \quad    i \in [1,s-1], \quad \bgamma \in   \bar{G}
\end{equation}
Consider $\gamma_s  $ with $  \bgamma \in   \bar{G}$. By  (\ref{4.0a}), we have
$$
   |\gamma_s  | =  | {\rm Nm} (\bgamma)  ( \gamma_1 \cdots \gamma_{s-1})^{-1} | \in   
 | \gamma_1 \cdots \gamma_{s-1} |^{-1} [C^s_M,\det \Gamma^\bot].
$$
Now using (\ref{4.0b}) and (\ref{4.0c}),  we obtain for sufficiently large $N$ that
\begin{equation}\nonumber
     \log_q |N_s\gamma_s  |  \leq n/s \log_q 2 +
       \log_q( c_0 \det \Gamma^\bot)     - m_1-...m_{s-1}   \leq  3/4n \log_q 2  ,
\end{equation}
and
\begin{equation}\nonumber
     \log_q |N_s \gamma_s  |  \geq   (n-1)/s \log_q 2  -\log_q c_0   + \log_q C_M^s - m_1-...m_{s-1} -s    \geq   n(\log_q 2)/(2s).
\end{equation}
 Therefore, we get for sufficiently large $N$ and for $  \bgamma \in   \bar{G}$ 
\begin{equation}\label{4.8}
  |\bgamma|<N/2, \; | {\rm Nm} (\bgamma)  | \leq  \det \Gamma^\bot ,\; |N_i \gamma_i  | \geq 2^{ln^4n}, \; i=1,...,s .
\end{equation}
So $\bar{G} \cup 2\bar{G} \subset G_1$ (see (\ref{3.21a})). 

Let  $  \bgamma \in   \bar{G}$.
Taking into account that $|N_i \gamma_i  | \geq  4$ $(i=1,...,s)$, we obtain
\begin{equation}\label{4.9}
  \int_0^1  \sin^2 (\pi\theta_i N_i \gamma_i)d\theta_i =1/2 - 1/2\int_0^1  \cos (2\pi\theta_i N_i \gamma_i) d\theta_i  \geq 1/4.
\end{equation}
Let $I = [1/6,1/3] \cup [2/3, 5/6] $. If $ \{\langle\bgamma,\bx\rangle    \}  \notin I $, then $|\cos(2\pi\{\langle\bgamma,\bx\rangle    \})| \geq 1/2$. 
Let $\{\langle\bgamma,\bx\rangle    \} \in I $. Then we take $2\bgamma$ instead of $\bgamma$. 
We see  that $|\cos(2\pi\{\langle2\bgamma,\bx\rangle    \})| \geq 1/2$, and
\begin{equation}\label{4.10}
    \max (\cos^2(2\pi\langle\bgamma(\bm),\bx\rangle  ),\cos^2(2\pi\langle2\bgamma(\bm),\bx\rangle  ) ) \geq 1/4.
\end{equation}
By   (\ref{4.0}) - (\ref{4.10}), we have
\begin{equation}\nonumber
  \check{ \psi} \geq   \sum_{\substack{\bgamma(\bm) \in \bar{G}\\ \bm \in \ZZ^{s-1}}} 
       \frac{(\det \Gamma)^{-2} }{(2\pi)^{2s} C_M^{2s}}  
       (1/4)^{s-1}  \max (\cos^2(2\pi\langle\bgamma(\bm),\bx\rangle  ),\cos^2(2\pi\langle2\bgamma(\bm),\bx\rangle  ) ) \geq w_4 \#\bar{G},
\end{equation}
with $w_4 =(\det \Gamma)^{-2} ((2\pi)^{2s} C_M^{2s})^{-1}  4^{-s}$. Applying  (\ref{4.3}),   we get the assertion of  Lemma~12. 
\\ \hfill \qed  \\

{\bf Lemma 13.} {\it There exist reals  $c_6,w_1>0$ such that for $N>c_6$}
\begin{equation}\label{4.50}
     {\bf E} [ (\cR (\btheta \cdot \bN  \cdot \KK^s+\bx, \Gamma) )^2 ]  \geq w_1 n^{s-1} .
\end{equation}

{\bf Proof.}  Applying  (\ref{3.10}) and  (\ref{3.23}), we have
\begin{equation}\label{4.51}
      {\bf E} [ | \cA(G^{(1)}, 0)|^2]  =  \ddot{S}_1 +  \ddot{S}_2,
\end{equation}
with
\begin{equation} \label{4.51a}
    \ddot{S}_1 =    \sum_{\bgamma^{(1)},\bgamma^{(2)}  \in G^{(1)}, \atop \bgamma^{(1)} =\pm\bgamma^{(2)}} 
    \ddot{\psi}(\bgamma^{(1)},-\bgamma^{(2)}), \qquad
        \ddot{S}_2     =    \sum_{\bgamma^{(1)} ,\bgamma^{(2)}  \in G^{(1)}, \atop \bgamma^{(1)} \neq \pm\bgamma^{(2)}} 
           \ddot{\psi}(\bgamma^{(1)},-\bgamma^{(2)}),
\end{equation}
where
\begin{equation}\nonumber
    \ddot{\psi}(\bgamma^{(1)},-\bgamma^{(2)}) =   
       \frac{(\det \Gamma)^{-2} \hat{\omega}(\tau\bgamma^{(1)}) \hat{\omega}(-\tau\bgamma^{(2)}) }
       {(2\pi)^{2s} {\rm Nm} (\bgamma^{(1)}) {\rm Nm} (-\bgamma^{(2)}) } e(\langle\bgamma^{(1)}-\bgamma^{(2)}\rangle,x)
        \boldsymbol{   \breve{\psi}}(\bgamma^{(1)}-\bgamma^{(2)}). 
\end{equation}
and
\begin{equation}\label{4.20}
   \boldsymbol{   \breve{\psi}}(\bgamma^{(1)},\bgamma^{(2)}) =           {\bf E} \Big[    \prod_{i=1}^s  \sin (\pi \theta_i N_i \gamma_i^{(1)})  
                \sin (-\pi\theta_i N_i\gamma_i^{(2)} ) \Big].
\end{equation}
We consider $\ddot{S}_1$.
 Bearing in mind (\ref{3}), (\ref{3.0}), that $|e(z)-1| =2|\sin(\pi z)| \leq 2\pi |z|$ and that $\omega(\bx)$ is supported inside 
the unit ball $B=\{\bx \; : \; |\bx|\leq 1\} $, we obtain for $\tau=1/N^2$ and $|\bgamma | \leq N$ that
\begin{equation}\label{4.52}
 |\hat{\omega} (\tau \bgamma) - 1| =         \Big|\int_{\RR^s}e(\langle \tau \bgamma,\bx\rangle) \omega(\bx)d\bx -1 \Big| =
    \Big|\int_{B }(e(\langle \bgamma,\bx\rangle/N^2)-1) \omega(\bx)d\bx \Big|
\end{equation}
\begin{equation}\nonumber
   \;\; \leq   \Big|\int_{B }|e(\langle \bgamma,\bx\rangle/N^2)-1| \omega(\bx)d\bx  \Big| 
 \leq 2 \pi/N \int_{B }\omega(\bx)d\bx = 2 \pi/N .
\end{equation}
By (\ref{4.20}), we see that
\begin{equation}\label{4.53}
  \boldsymbol{   \breve{\psi}}(\bgamma^{(1)},-\bgamma^{(2)})\big( {\rm Nm} (\bgamma^{(1)}) {\rm Nm} (-\bgamma^{(2)})  \big)^{-1}   
   =  \boldsymbol{   \breve{\psi}}(\bgamma^{(1)},\bgamma^{(2)}) \big( {\rm Nm} (\bgamma^{(1)}) {\rm Nm} (\bgamma^{(2)})  \big)^{-1}  .
\end{equation}
Taking into account that $1 + \cos(2z) =2\cos^2(z)$, we get from   (\ref{4.0}), (\ref{4.51a}), (\ref{4.52}) and (\ref{4.53})  that
\begin{equation}\label{4.52b}
   \ddot{S}_1 = \check{ \psi} +O\Big(1/N  \sum_{\bgamma \in G^{(1)}} 
       \frac{1}{ |{\rm Nm} (\bgamma)|^2}  \Big).
\end{equation}
By (\ref{3.300}), (\ref{3.37}), (\ref{4.52})  and Lemma 12, we have for sufficiently large $N$ 
\begin{equation}\label{4.531}
   \ddot{S}_1 \geq 0.5 w_3n^{s-1}.
\end{equation}
Now we consider $\ddot{S}_2$.
We see from (\ref{4.20}) and (\ref{3.262})  
\begin{equation}\nonumber
   |\boldsymbol{   \breve{\psi}}(\bgamma^{(1)},\bgamma^{(2)})|  = 
   \Big|2^{-2s} \sum_{\kappa_1^{(1)},...,\kappa_s^{(2)} \in \{-1, 1\} } \kappa_1^{(1)}\kappa_2^{(1)} \cdots \kappa_s^{(2)}
             {\bf E} \Big[ e\Big(  \sum_{1 \leq i \leq s, j=1,2}    \theta_i N_i \kappa_i^{(j)}\gamma_i^{(j)}/2  \Big)\Big] \Big|  
\end{equation} 
\begin{equation}\nonumber
    \leq  
  2^{-s} \sum_{\kappa_1,...,\kappa_s \in \{-1, 1\}} 
             \Big|{\bf E} \Big[ e\Big(  \sum_{i=1}^s   \theta_i N_i (\gamma_i^{(1)} +  \kappa_i \gamma_i^{(2)})/2 \Big)\Big] \Big| 
\end{equation} 
\begin{equation}\label{4.22}
    \leq  
  2^{-s} \sum_{\kappa_1,...,\kappa_s \in \{-1, 1\}} 
            \prod_{i=1}^s  \min\Big(1, \frac{1}{\pi | N_i (\gamma_i^{(1)} +  \kappa_i \gamma_i^{(2)})|} \Big). 
\end{equation} 
Applying Lemma 3, we get that $\psi(\bgamma^{(1)},\bgamma^{(2)}) =O(n^{-20s})$.
By (\ref{3.305}), (\ref{3.38}) and  (\ref{4.51a}), we derive that
\begin{equation}\nonumber   
   \ddot{S}_2 =    O(n^{-20s} S_2(G^{(1)}, G^{(1)})) =O(n^{-2s})  .
\end{equation}
From (\ref{4.51}) and (\ref{4.531}),   we have for sufficiently large $N$ that
\begin{equation}\label{4.532}
  {\bf E} [ | \cA(G^{(1)},0)|^2]   \geq 0.25 w_3n^{s-1}.
\end{equation}
By the triangle inequality, we obtain
\begin{equation}\nonumber
        \big( {\bf E} [ (\cR (\btheta \cdot \bN  \cdot \KK^s+\bx, \Gamma) )^2 ]  \big)^{1/2}  \geq    \big(  {\bf E} [ | \cA(G^{(1)}, 0)|^2] \big)^{1/2}
\end{equation}
\begin{equation}\nonumber
             -    \big(  {\bf E} [| \cR (\btheta \cdot \bN  \cdot \KK^s+\bx, \Gamma) -  \cA(G^{(1)},0) |^2]\big)^{1/2}.
\end{equation}
Using (\ref{4.532}) and Lemma 11, we get the assertion of Lemma 13. \qed  \\

\setcounter{equation}{0}
\renewcommand{\theequation}{2.6.\arabic{equation}}

{\bf 2.6. Four moments estimates  for $ \cA(\dot{G}_i)$.} In this subsection, we will prove that 
\begin{equation}\label{6.2}
      {\bf E} [ |  \cA(\dot{G}_i)|^4] = O( i^{2(s-2)}n^{8/9(s-1)})  
\end{equation}
and 
\begin{equation}\nonumber
    {\bf E} \Big[ \Big( \sum_{i \in [1,n^{5/9}]}   (\cA^2(\dot{G}_i) -  {\bf E}  [  \cA^2(\dot{G}_i)] )^2   \Big)\Big] =O(n^{2(s-1)-2/9}) .
\end{equation}\\
We need these estimates to apply the martingale CLT in the next section.
Let
\begin{equation} \nonumber
   \delta(\fT) =   \begin{cases}
    1,  & \; {\rm if}  \;  \fT  \;{\rm is \;true},\\
    0, &{\rm otherwise}.
  \end{cases}
\end{equation}
 
{\bf Lemma 14.} {\it The assertion (\ref{6.2}) is true. } \\

{\bf Proof.} 
Using  the following simple inequality
\begin{equation} \nonumber
      \Big| \sum_{1 \leq i \leq 2^s} a_i\Big|^4  \leq   ( 2^s \max_{1 \leq i \leq 2^s} |a_i|)^4 \leq   2^{4s} \sum_{1 \leq i \leq 2^s} |a_i|^4,
\end{equation}
 we obtain from (\ref{3.232})
\begin{equation}\label{6.3}
      {\bf E} [ |  \cA(\dot{G}_i,\tau)|^4]  \leq  2^{4s} \sum_{\kappa_1,...,\kappa_s \in \{-1, 1\}} |  \cB(\dot{G}_i,\bkappa,\tau)|^4.
\end{equation}
Applying  (\ref{3.24b}) and Lemma 6, we get
\begin{equation} \nonumber
      {\bf E} [   |\cB((\dot{G}_i,\bkappa)|^4)] =    \prod_{1 \leq j \leq 4}
       \sum_{\bgamma^{(j)} \in  \dot{G}_i}
       | h(\bgamma^{(j)}) | ( \delta(\hat{\bgamma} =\b0) + (1-  \delta(\hat{\bgamma} =\b0))N^{-1}O(1)) , 
\end{equation}
where $ \hat{\bgamma} =  \bgamma^{(1)} - \bgamma^{(2)}            +  \bgamma^{(3)}       - \bgamma^{(4)}$.
 From (\ref{3.24a}), (\ref{3.305}), (\ref{3.38}) and Lemma 2,  we derive that
\begin{equation}\label{6.4a}
      {\bf E} [   |\cB((\dot{G}_i,\bkappa,\tau)|^4)] =    V_{1} + V_{2}   + O(n^{8s}/N),    
\end{equation}
where 
\begin{equation}\label{6.4a1}
     V_{k} =    \prod_{1 \leq j \leq 4}
       \sum_{\bgamma^{(j)} \in  \dot{G}_i}
             \frac{(\det \Gamma)^{-1} |\hat{\omega}(\tau\bgamma^{(j)})|}{(2\pi)^s |{\rm Nm} (\bgamma^{(j)})|}  
             \delta(\hat{\bgamma} =\b0)
             \delta_k(\bar{\bgamma}), \quad k=1,2,
\end{equation}
with  $\bar{\bgamma} = (\bgamma^{(1)},\bgamma^{(2)},\bgamma^{(3)},\bgamma^{(4)})$ and
\begin{equation}\label{6.4c}
         \delta_1(\bar{\bgamma}) =\delta( \nexists j,l \in [1,4], \; j \neq l \; | \; \bgamma^{(j)} =(-1)^{l-j+1}\bgamma^{(l)}  ), 
           \qquad
             \delta_2(\bar{\bgamma}) = 1-\delta_1(\bar{\bgamma}) .
\end{equation}
By  (\ref{2.82}) and (\ref{2.9}), we have that $\bgamma^{(j)} = \tilde{\bgamma}^{(j)} \cdot \sigma(\eta^{(j)})$, with
 $\tilde{\bgamma}^{(j)}=(\tilde{\gamma}_1^{(j)},..., \tilde{\gamma}_s^{(j)}) $,   $\tilde{\bgamma}^{(j)} \in \FF_n$  and $\eta^{(j)} \in \dot{\fU}({\bgamma}^{(j)},a,b)$
   $(j=1,...,4)$.  Using  (\ref{2.9}) and (\ref{3.56}), we obtain that $a=i n^{4/9} $ and $ b=n^{4/9} - n^{2/9}$.\\
Hence
\begin{equation}\label{6.5}
        V_{k} =       O\Big(    \prod_{1 \leq j \leq 4}  \sum_{\tilde{\bgamma}^{(j)} \in \FF_n} 
        \sum_{\eta^{(j)} \in \dot{\fU}(\tilde{\bgamma}^{(j)},a,b)} 
       \frac{1 }{ |{\rm Nm} ( \tilde{\bgamma}^{(j)}) |}  \delta_k(\bar{\bgamma})
  \Big) 
\end{equation}
$$
  \times \delta(   \tilde{\bgamma}^{(1)}  \cdot \sigma(\eta^{(1)}) - \tilde{\bgamma}^{(2)}\cdot \sigma(\eta^{(2)})
            +  \tilde{\bgamma}^{(3)} \cdot \sigma(\eta^{(3)})       - \tilde{\bgamma}^{(4)}\cdot \sigma(\eta^{(4)})  =\b0).             
$$
 It is easy to see that
$$
     \tilde{\bgamma}^{(1)} \cdot \sigma(\eta^{(1)}) - \tilde{\bgamma}^{(2)}\cdot \sigma(\eta^{(2)})
            +  \tilde{\bgamma}^{(3)} \cdot \sigma(\eta^{(3)})       - \tilde{\bgamma}^{(4)}\cdot \sigma(\eta^{(4)})  =\b0
$$
if and only if
\begin{equation} \label{6.5a}              
    (\tilde{\gamma}_1^{(1)}\sigma_1(\eta^{(1)})-
               \tilde{\gamma}_1^{(2)}\sigma_1(\eta^{(2)})
                  +  \tilde{\gamma}_1^{(3)} \sigma_1(\eta^{(3)}) )/ \tilde{\gamma}_1^{(4)}\sigma_1(\eta^{(4)})      =1.
\end{equation}
First we consider $V_{1}$. 
    We fix  $  \tilde{\bgamma}^{(1)},\tilde{\bgamma}^{(2)},  \tilde{\bgamma}^{(3)}, \tilde{\bgamma}^{(4)}$ and $\eta^{(4)} $.
              From (\ref{2.82}) and (\ref{6.4c}), we get that there is no degenerate solutions $(\eta^{(1)},\eta^{(2)}, \eta^{(3)})$ of the equation 
              (\ref{6.5a}).
Applying Theorem A, we have that the number of non-degenerate solutions $(\eta^{(1)},\eta^{(2)}, \eta^{(3)})$ of 
              (\ref{6.5a}) is finite. Hence
\begin{equation}\label{6.4f}
     V_1 =    O\Big( 
      \sum_{\tilde{\bgamma}^{(j)} \in \FF_n, 1 \leq j \leq 4} 
        \sum_{\eta^{(4)} \in \dot{\fU}(\tilde{\bgamma}^{(4)},a,b)} 
       \frac{1 }{ |{\rm Nm} (\tilde{\bgamma}^{(1)})| \cdots |{\rm Nm} (\tilde{\bgamma}^{(4)})| }       \Big).
\end{equation}
By    (\ref{2.10}) and  (\ref{2.10a}), we derive
\begin{equation}\label{6.4g}
     V_{1} =   O((\ln \; n)^4    b(a+b)^{s-2}) =O(i^{s-2}n^{4/9(s-1)}(\ln \; n)^4  ) .
\end{equation}
Now we consider $V_2$.  Let $ \bgamma^{(j_0)} =(-1)^{l_0-j_0+1} \bgamma^{(l_0)}  $.
Bearing in mind that $\hat{\bgamma} =  \bgamma^{(1)} - \bgamma^{(2)}
            +  \bgamma^{(3)}       - \bgamma^{(4)} =\b0 $, we obtain that $ \bgamma^{(j_1)} =(-1)^{l_1-j_1+1} \bgamma^{(l_1)} $
with $\{j_1,l_1\} = \{1,2,3,4\}   \setminus \{ j_0,l_0\}$. 
Hence, from (\ref{6.4a1}), we get
\begin{equation}\nonumber
        V_{2} =       O\Big( \Big(  \sum_{\tilde{\bgamma}^{(1)} \in \FF_n} 
        \sum_{\eta^{(1)} \in \dot{\fU}(\tilde{\bgamma}^{(1)},a,b)} 
       \frac{1 }{ {\rm Nm}^2 (\tilde{\bgamma}^{(1)})}
  \Big)^2 \Big)  .
\end{equation}
By  (\ref{2.10}) and  (\ref{2.10a}), we have
\begin{equation}\nonumber
     V_{2} =   O(    b^2(a+b)^{2(s-2)}) =O( i^{2(s-2)}n^{8/9(s-1)}) .
\end{equation}
Using   (\ref{6.3}),  (\ref{6.4a}) and  (\ref{6.4g}), we obtain  (\ref{6.2}) and the assetetion of Lemma 14. \qed \\

Let  $\quad \bar{\bgamma} = (\bgamma^{(1)},\bgamma^{(2)},\bgamma^{(3)},\bgamma^{(4)})$, and let
\begin{equation}\nonumber
    \dot{\bgamma} =  \dot{\bgamma} =(\dot{\gamma}_1,...,\dot{\gamma}_s ) =\bkappa^{(1)} \cdot \bgamma^{(1)} + \cdots + \bkappa^{(4)} \cdot \bgamma^{(4)}  ,
       \quad
       \ddot{\bgamma} =  \bkappa^{(3)} \cdot \bgamma^{(1)} +  \bkappa^{(4)} \cdot \bgamma^{(2)},
\end{equation}
\begin{equation}\label{6.00}
       \ddot{\bgamma} =  \ddot{\bgamma} =(\ddot{\gamma}_1,...,\ddot{\gamma}_s ) = \bkappa^{(3)} \cdot \bgamma^{(1)} +  \bkappa^{(4)} \cdot \bgamma^{(2)},
\end{equation}

\begin{align} \nonumber
   & \dot{\delta}_1(\bar{\bgamma}) = \delta(  \dot{\bgamma}= \b0, \;  \ddot{\bgamma}= \b0), 
   \\
     \nonumber
   & \dot{\delta}_2(\bar{\bgamma}) = \delta(  \dot{\bgamma} =\b0, \; \ddot{\bgamma}  \neq 0), 
    \\
 \nonumber
   & \dot{\delta}_3(\bar{\bgamma}) = \delta(  \dot{\bgamma} \neq \b0, \; \ddot{\bgamma}= \b0),
   \\
\nonumber
   & \dot{\delta}_4(\bar{\bgamma}) = \delta(  \dot{\bgamma} \neq \b0, \; \ddot{\bgamma} \neq \b0, \;\;\qquad \qquad\qquad \quad \exists \nu \in [1,s]  \;\;\; \dot{\bgamma}_{\nu}  = 0 \;\;\;
    \quad  {\rm and } \;\;\;\quad \ddot{\bgamma}_{\nu}  = 0  ),  \\
\nonumber
    & \dot{\delta}_5(\bar{\bgamma}) = \delta( \dot{\bgamma} \neq \b0, \; \ddot{\bgamma} \neq \b0,\;\; \quad \exists \nu \in [1,s]  \;\;
       \dot{\bgamma}_{\nu}  = 0 \;\; 
         {\rm and } \;\; \nexists \nu \in [1,s] \;\;\; \dot{\bgamma}_{\nu} = 0,\; \ddot{\bgamma}_{\nu} = 0),  \\
\nonumber
   & \dot{\delta}_6(\bar{\bgamma}) = \delta(  \dot{\bgamma} \neq \b0, \;  \ddot{\bgamma} \neq \b0, \;\; \qquad\qquad\qquad \quad \forall \nu \in [1,s]  \;\; \dot{\bgamma}_{\nu}  \neq 0
    \quad \;\;\;  {\rm and } \;\;\;\;\quad \ddot{\bgamma}_{\nu} \neq 0),  \\
\nonumber
   & \dot{\delta}_7(\bar{\bgamma}) = \delta(  \dot{\bgamma} \neq \b0, \; \ddot{\bgamma} \neq \b0,\;\; \quad \forall \nu \in [1,s]  \;\;\; \dot{\bgamma}_{\nu}  \neq 0 
     \;\;\; \quad         {\rm and } \;\;\;\quad \exists \nu \in [1,s] \;\;\; \;\; \;\ddot{\bgamma}_{\nu} = 0).
\end{align}
It is easy to verify that
\begin{equation}\label{6.07}
     \sum_{1 \leq k \leq 7}   \dot{\delta}_k(\bar{\bgamma}) =1.
\end{equation} \\

 {\bf Lemma 15.} {\it Let $l \geq 2,\; i,j=1,2,...$, and let
\begin{equation}\label{6.08}
   \dot{H}_{i,j}(l)= \sum_{\bkappa^{(1)},...,\bkappa^{(4)} \in \{-1, 1\}^s}   \sum_{\bgamma^{(1)},\bgamma^{(2)} \in \dot{G}_i} 
    \sum_{\bgamma^{(3)},\bgamma^{(4)} \in \dot{G}_j }        h(\bgamma^{(1)})\cdots h(\bgamma^{(4)})    
\end{equation} 
\begin{equation}\nonumber
        \times  
             \dot{\delta}_l(\bar{\bgamma})  {\bf E} \Big[ e\Big(\sum_{k \in [1,s] } \sum_{\nu  \in [1,2]}\theta_k \kappa_k^{(\nu)} \gamma_k^{(\nu)}  N_k /2\Big)\Big] 
           {\bf E} \Big[e\Big(\sum_{k \in [1,s] } \sum_{\nu  \in [3,4]} \theta_k \kappa_k^{(\nu)} \gamma_k^{(\nu)}  N_k /2 \Big) \Big]   .
\end{equation} 
Then }
\begin{equation}\label{6.081a}
       \dot{H}_{i,j}(l)  =O(n^{-10s}).
\end{equation} \\

{\bf Proof.} Applying  (\ref{3.24a}) and Lemma 6, we get
\begin{equation} \nonumber
   \dot{H}_{i,j}(l)  =O\Big(  \sum_{\bkappa^{(1)},...,\bkappa^{(4)} \in \{-1, 1\}^s}   \sum_{\bgamma^{(1)},\bgamma^{(2)} \in \dot{G}_i} 
    \sum_{\bgamma^{(3)},\bgamma^{(4)} \in \dot{G}_j }      |{\rm Nm} (\bgamma^{(1)}...\bgamma^{(4)})|^{-1}     \dot{\delta}_l(\bar{\bgamma})
\end{equation}
\begin{equation} \nonumber
    \times \prod_{k \in [1,s]} \min(1, \frac{1}{N_k|\dot{\gamma}_k  - \ddot{\gamma}_k|})\min(1, \frac{1}{N_k|\ddot{\gamma}_i|})\Big).
\end{equation}
From (\ref{6.00}), we have that for $l \geq 2$ there exists $k_0 \in [1,s]$ such that $\max(|\dot{\gamma}_{k_0}  - \ddot{\gamma}_{k_0}|, | \ddot{\gamma}_{k_0}|) >0$. Using Lemma 3, we derive that $N_i \max(|\dot{\gamma}_{k_0}  - \ddot{\gamma}_{k_0}|, | \ddot{\gamma}_{k_0}|) >\dot{c} n^{20s}$.
%c_1 2^{(\ln n)^4}$.
 Thus
\begin{equation}\label{6.073}
   \dot{H}_{i,j}(l)  =O\Big(   n^{-20s} \Big(   \sum_{\bgamma \in G^{(1)}} 
     ({\rm Nm} (\bgamma))^{-1}    \Big)^4 \Big).
\end{equation}
By    (\ref{2.90}) and  (\ref{3.21a})  $G^{(1)} = \GG(0,2n)$.  Similarly to (\ref{6.4f}) and  (\ref{6.4g}), we 
obtain from Lemma 2  that
\begin{equation}\label{6.073b}
   \dot{H}_{i,j}(l)  =O(  n^{-20s} n^{8s}) =O(n^{-10s}).
\end{equation}
Hence,  Lemma 15  is proved. \qed \\

{\bf Lemma 16.} {\it Let  $l \geq 2, \; i <j$, and
\begin{equation}\label{6.071}
   \ddot{H}_{i,j}(l)= \sum_{\bkappa^{(1)},...,\bkappa^{(4)} \in \{-1, 1\}^s}   \sum_{\bgamma^{(1)},\bgamma^{(2)} \in \dot{G}_i} 
    \sum_{\bgamma^{(3)},\bgamma^{(4)} \in G_j }        h(\bgamma^{(1)})\cdots h(\bgamma^{(4)})    
\end{equation}
\begin{equation}\nonumber
        \times     \dot{\delta}_l(\bar{\bgamma})
            {\bf E} \Big[e\Big(\sum_{k \in [1,s] } \sum_{\nu  \in [1,4]}  \theta_k \kappa_k^{(\nu)} \gamma_k^{(\nu)}  N_k /2 \Big) \Big]  . 
\end{equation}
Then}
\begin{equation}\label{6.072}
      \dot{H}_{i,j}(l)  =O(j^{s-2}n^{4/9(s-1)+ 2/45} )   .
\end{equation}  \\

 {\bf Proof.}     Applying  (\ref{3.24a}) and Lemma 6, we get
\begin{equation}\label{6.073m}
   \ddot{H}_{i,j}(l)  =O\Big(  \sum_{\bkappa^{(1)},...,\bkappa^{(4)} \in \{-1, 1\}^s}   \sum_{\bgamma^{(1)},\bgamma^{(2)} \in \dot{G}_i} 
    \sum_{\bgamma^{(3)},\bgamma^{(4)} \in \dot{G}_j }      |{\rm Nm} (\bgamma^{(1)}...\bgamma^{(4)})|^{-1}    
\end{equation}
\begin{equation}\nonumber
    \times  \dot{\delta}_l(\bar{\bgamma})   \prod_{\nu \in [1,s]} \min(1, \frac{1}{N_{\nu}|\dot{\gamma}_{\nu}|})\Big).
\end{equation}
We will prove Lemma 16 separately for each $l \in [2,7]$:

 { \em Case} $l \in \{2,5\}$.  We will consider the case $l=2$. The proof for the case $l=5$ is similar.
  By (\ref{6.00}) and  (\ref{6.073m}),  we have
\begin{equation}\nonumber
   \dot{H}_{i,j}(2)  =O(    \sum_{1 \leq \nu \leq s}\tilde{H}_{i,j,\nu}).
\end{equation}
with
\begin{equation}\nonumber
 \tilde{H}_{i,j,\nu} = \sum_{\bkappa^{(1)},...,\bkappa^{(4)} \in \{-1, 1\}^s}   \sum_{\bgamma^{(1)},\bgamma^{(2)} \in \dot{G}_i} 
    \sum_{\bgamma^{(3)},\bgamma^{(4)} \in G_j }      |{\rm Nm} (\bgamma^{(1)}...\bgamma^{(4)})|^{-1}    \delta( \dot{\gamma}_{\nu}  =0,
    \ddot{\gamma}_{\nu}  \neq 0).
\end{equation}
Let $\bgamma^{(j)} = \tilde{\bgamma}^{(j)} \cdot \sigma(\eta^{(j)})$, with
   $\tilde{\bgamma}^{(j)} \in \FF_n$  and $\eta^{(j)} \in \dot{\fU}({\bgamma}^{(j)},a_j,b_j)$
   $(j=1,...,4)$,  $a_1=a_2=i n^{4/9}, \; a_3=a_4=j n^{4/9}  $ and $ b_1=...=b_4=n^{4/9} - n^{2/9}$.
    We fix  $  \tilde{\bgamma}^{(1)},\tilde{\bgamma}^{(2)},  \tilde{\bgamma}^{(3)}, \tilde{\bgamma}^{(4)}$ and $\eta^{(4)} $.
 Bearing in mind that $\ddot{\gamma_{\nu}} \neq 0$ and $i <j$,  we obtain that there is no degenerate solutions $(\eta^{(1)},\eta^{(2)}, \eta^{(3)})$ (see (\ref{6.0})) of the equation 
$$              
          \dot{\gamma}_{\nu}  =    \sum_{1 \leq k \leq 4} \kappa_{\nu}^{(k)} \tilde{\gamma}_{\nu}^{(k)} \sigma_{\nu}(\eta^{(k)}) =0 .
$$ 
Similarly to (\ref{6.4f}) and (\ref{6.4g}),  we derive from Theorem A,  (\ref{2.10}) and  (\ref{2.10a}) that 
\begin{equation}\nonumber
     \tilde{H}_{i,j,\nu} =    O\Big(       \sum_{\tilde{\bgamma}^{(j)} \in \FF_n, 1 \leq j \leq 4} 
        \sum_{\bet^{(4)} \in \dot{\fU}(\tilde{\bgamma}^{(4)},a_4,b_4)} 
       \frac{1 }{ {\rm Nm} (\tilde{\bgamma}^{(j)}) }       \Big)  =O(j^{s-2}n^{4/9(s-1)}(\ln \; n)^4  ).
\end{equation}
Hence, the assertion (\ref{6.072}) is proved.

 { \em Case} $l\in \{3,7\}$. We have from (\ref{6.00}) for  both cases $l=3$ and $l=7$ that there exists $\nu \in [1,s]$ such that $   \dot{\gamma}_{\nu} \neq 0$, $\ddot{\gamma}_{\nu}=0$
 and   $\dot{\gamma}_{\nu}  =     \kappa_{\nu}^{(1)} \gamma_{\nu}^{(1)} +  \kappa_{\nu}^{(2)} \gamma_{\nu}^{(2)}$.
 Applying Lemma 3, we get  $|N_{\nu}\dot{\gamma}_{\nu}| \geq \dot{c} n^{20s} $. Now using (\ref{6.073}),(\ref{6.073b}) and (\ref{6.073m}), we obtain
 (\ref{6.072}).

{ \em Case} $l=4$. By (\ref{6.00}), we have that there exist $\mu,\nu \in [1,s]$ with $   \ddot{\gamma}_{\nu} = 0$,
 $   \dot{\gamma}_{\nu} = 0$, and  $\dot{\gamma}_{\mu} \neq  0$.
It is easy to derive  that  $ \bgamma^{(1)} = \pm  \bgamma^{(2)}$, $ \bgamma^{(3)} = \pm  \bgamma^{(4)}$  and
  $\dot{\gamma}_{\mu} = \tilde{\kappa}_1 \gamma_{\mu}^{(1)} + \tilde{\kappa}_2 \gamma_{\mu}^{(3)}$  with $\tilde{\kappa}_i \in \{-2,0,2\}$,
 $i=1,3$. Hence
$$  
  |\dot{\gamma}_{\mu}| = 2|\gamma_{\mu}^{(1)}| \quad {\rm or}  \quad    |\dot{\gamma}_{\mu}| = 2|\gamma_{\mu}^{(3)}| \quad  {\rm or} 
  \quad |\dot{\gamma}_{\mu}| = 2 |  \gamma^{(1)}_{\mu} \pm \gamma^{(3)}_{\mu}| \neq 0.
$$
Applying  (\ref{3.21a}) and Lemma 3, we get $ |N_{\mu}\dot{\gamma}_{\mu}| \geq \dot{c}n^{20s} $ for sufficiently large $N$. By
 (\ref{6.073}), (\ref{6.073b}) and  (\ref{6.073m}), we obtain
 (\ref{6.072}).

{ \em Case}  $l=6$. By (\ref{3.56}), we have that there exists $\nu \in [1,s]$ such that  $\gamma^{(4)}_{\nu}  \geq 2^{j n^{4/9}} $.
Using Lemma 3, we obtain  for sufficiently large $N$ that
\begin{equation}\nonumber
       |\ddot{\gamma}_{\nu}| =       | \kappa^{(3)}_{\nu} \gamma^{(3)}_{\nu}   +  \kappa^{(4)}_{\nu}  \gamma^{(4)}_{\nu} |
             \geq |\gamma_{\nu}^{(4)}|   \exp( -\ddot{c}  (\ln n)^3) \geq 2^{j n^{4/9}}  \exp( -\ddot{c}  (\ln n)^3)|     
\end{equation} 
\begin{equation}\nonumber
   \geq 2^{(i+1) n^{4/9} - n^{2/9}+2}  
              \geq
             2 |\kappa^{(1)}_{\nu} \gamma^{(1)}_{\nu}   +  \kappa^{(2)}_{\nu}  \gamma^{(2)}_{\nu}|.
\end{equation} 
Hence, we get  for sufficiently large $N$ that
\begin{equation}\nonumber
       N_{\nu}   |\dot{\gamma}_{\nu}| =     N_{\nu}  | \kappa^{(1)}_{\nu} \gamma^{(1)}_{\nu}   +  \kappa^{(2)}_{\nu}  \gamma^{(2)}_{\nu} + 
               \kappa^{(3)}_{\nu} \gamma^{(3)}_{\nu}   +  \kappa^{(4)}_{\nu}  \gamma^{(4)}_{\nu} | \geq N_{\nu}  |\ddot{\gamma}_{\nu}|/2
             \geq  n^{20s}.  
\end{equation} 
 Now from (\ref{6.073}), (\ref{6.073b})  and  (\ref{6.073m}), we obtain
 (\ref{6.072}). 
Thus, Lemma 16 is proved. \hfill \qed 
\\

Let
\begin{equation}\label{6.27a}
   H_{i,j}=  \sum_{1 \leq l \leq 7} ( \ddot{H}_{i,j}(l) -   \dot{H}_{i,j}(l)  ).   
\end{equation}
By    (\ref{6.07}), (\ref{6.08}) and  (\ref{6.071}), we get
\begin{equation}\label{6.25}
   H_{i,j}= \sum_{\bkappa^{(1)},...,\bkappa^{(4)} \in \{-1, 1\}^s}   \sum_{\bgamma^{(1)},\bgamma^{(2)} \in \dot{G}_i} 
    \sum_{\bgamma^{(3)},\bgamma^{(4)} \in G_j }        h(\bgamma^{(1)})\cdots h(\bgamma^{(4)})    
\end{equation}
\begin{equation}\nonumber
        \times  
            \Big({\bf E} \Big[e(\sum_{k \in [1,s] } \sum_{l  \in [1,4]}  \phi_{k,l})\Big] - 
            {\bf E} \Big[e(\sum_{k \in [1,s] } \sum_{l  \in [1,2]} \phi_{k,l})\Big]
             {\bf E} \Big[e(\sum_{k \in [1,s] } \sum_{l  \in [3,4]} \phi_{k,l})\Big] \Big)     
\end{equation}
with $\phi_{k,l} = \theta_k \kappa_k^{(l)} \gamma_k^{(l)}  N_k /2$. \\

 {\bf Lemma 17.} {\it With notations as above, we have}
\begin{equation}\nonumber
   \varkappa:=   {\bf E} \Big[ \Big( \sum_{i \in [1,n^{5/9}]}   (\cA^2(\dot{G}_i) -  {\bf E}  [  \cA^2(\dot{G}_i)] )  \Big)^2 \Big]
    =O(   n^{2(s-1)-2/5} ) . 
\end{equation}\\
{\bf Proof.}  Let
\begin{equation}\nonumber
    \varkappa_{i,j} = {\bf E}\big[ \big(  \cA^2(\dot{G}_i)
     -  {\bf E}   [\cA^2(\dot{G}_i)] \big)  \times \big(   \cA^2(\dot{G}_j) -  {\bf E}  [  \cA^2(\dot{G}_j)]\big) \big] . 
\end{equation}
It is easy to see that
\begin{equation}\nonumber
    \varkappa_{i,j} = {\bf E}\big[  \cA^2(\dot{G}_i) \cA^2(\dot{G}_j)] 
     -  {\bf E}   [\cA^2(\dot{G}_i)] {\bf E}  [  \cA^2(\dot{G}_j)],  
\end{equation}
and
\begin{equation}\label{6.25a1}
   \varkappa  \leq   \dot{\varkappa} + \ddot{\varkappa}, \; {\rm with} \;
    \dot{\varkappa} =  \sum_{i \in [0,n^{5/9}]}  \varkappa_{i,i}    \; {\rm and} \;  
       \ddot{\varkappa} =
           2\sum_{i,j \in [1,n^{5/9}], \; i < j}     |\varkappa_{i,j}| .
\end{equation}
By Lemma 14, we obtain
\begin{equation} \nonumber
     \varkappa_{i,i}  \leq {\bf E} [ |  \cA(\dot{G}_i)|^4] = O( i^{2(s-2)}n^{8/9(s-1)}),  
\end{equation}
and
\begin{equation}\label{6.25a2}
     \dot{\varkappa} =O\Big( \sum_{i \in [0,n^{5/9}]}  i^{2(s-2)}n^{8/9(s-1)}) = O(   n^{2(s-1)-5/9} \Big).
\end{equation}
Using  (\ref{3.24c}) and (\ref{6.25}),     we get
\begin{equation}\label{6.25a}
    \varkappa_{i,j} = H_{i,j}.   
\end{equation}
From  (\ref{6.00}), (\ref{6.08}) and (\ref{6.071}),   we derive 
$$
            \ddot{H}_{i,j}(1) - \dot{H}_{i,j}(1)    =0.
$$
By Lemma 15 and Lemma 16,    we have
$$
  \dot{H}_{i,j}(l) =O(n^{-10s})  \quad {\rm and} \quad
    \ddot{H}_{i,j}(l) =O(j^{s-2}n^{4/9(s-1)+2/45} ), \quad l=2,3,...,7, \; i<j.
$$ 
Applying  (\ref{6.27a}), we obtain  
$
 H_{i,j}=  O(j^{s-2}n^{4/9(s-1)}).
$ 
Now from (\ref{6.25a1}) and (\ref{6.25a}), we get
\begin{equation}\nonumber
   \ddot{\varkappa} =O\Big( \sum_{j \in [1,n^{5/9}]}  j^{s-1}n^{4/9(s-1)+2/45}   \Big)=O(  n^{5s/9 +(s-1)4/9+2/45}  ) = O(   n^{s-2/5} ) .
\end{equation}  
By (\ref{6.25a1})   and  (\ref{6.25a2}),  Lemma 17 is proved. \qed \\

\setcounter{equation}{0}
\renewcommand{\theequation}{2.7.\arabic{equation}}
{\bf 2.7. Martingale approximation.}

  Denote by $\dot{\cF}(l)$ the sigma field on $[0,1)^s$ generated by 
  $\{ [\frac{k_1}{2^l},\frac{k_1+1}{2^l}) \times ... \times [\frac{k_s}{2^l},\frac{k_s+1}{2^l} ) \; : \;
   k_1,...,k_s =0,..., 2^l -1  \}$. Let $l(0)=0, l(i)=(i+1)[n^{4/9}] + [n/s - n^{1/9}],$
\begin{equation}\label{5.1}
 \cF_i=\dot{\cF}(l(i))   \qquad   {\rm and}  \qquad  
      \xi_i = {\bf E} [ \cA(\dot{G}_i ) \; | \; \cF_{i}] - {\bf E} [ \cA(\dot{G}_i ) \; | \; \cF_{i-1}], \quad i=1,2,...
\end{equation}   
Then $(\xi_i)_{i \geq 1}$ is the martingale difference array satisfying ${\bf E} [ \xi_i  |  \cF_{i-1}]=0, \; i=1,2,...$\\  

{\bf Lemma 18.} {\it With notations as above  
\begin{equation}\label{5.2}
      {\bf E} [ \cA(\dot{G}_i ) \; | \; \cF_{i-1}]  =O(n^{-10s}),   \quad \quad  \quad  \quad  \cA(\dot{G}_i ) - \xi_i = O(n^{-10s}) ,
\end{equation}
\begin{equation}\label{5.21}
  \cA(\dot{G}_i )^2 - \xi_i^2 = O(n^{-8s})      \quad \quad {\it and}  \quad  \quad  |\xi_i|^4 \leq  8|\cA(\dot{G}_i )|^4+ O(n^{-6s}) .
\end{equation}}
{\bf Proof.} It is easy to see that
\begin{equation}\nonumber
     \big|    2^l \int_{k/ 2^l}^{(k+1)/ 2^l} \sin (\lambda  \theta)  d \theta  \big|  \;
    = \; 2^l/\lambda |    \cos (\lambda  (k+1)/ 2^l) -  \cos (\lambda  k/ 2^l)|  
  \leq  2^{l+1}/\lambda ,\quad {\rm with} \quad \lambda >0.
\end{equation}
Hence, we obtain for $|\gamma_j| \geq 2^{i[n^{4/9}]}$ and $|N_j | \geq c_0^{-1}2^{(n-1)/s}$ that
\begin{equation}\label{5.4}
     \Big|    2^{l_{i-1}} \int_{k/ 2^{l_{i-1}}}^{(k+1)/ 2^{l_{i-1}}} \sin (N_j \gamma_j \theta  )  d \theta  \Big|  
  \leq  c_0 2^{-   [n^{1/9}] +5}. 
\end{equation}
Bearing in mind that
\begin{equation}\nonumber
  {\bf E} [  \phi_1(\theta_1 ) \cdot... \cdot  \phi_s(\theta_s ) \; | \; \cF_{i-1}]  = 
   \prod_{j=1}^s    2^{l_i} \int_{k_j/ 2^{l_{i-1}}}^{(k_j+1)/ 2^{l_{i-1}}} \phi_j(\theta_j )    d \theta_j 
\end{equation}
on  $ [\frac{k_1}{2^{l_{i-1}}},\frac{k_1+1}{2^{l_{i-1}}})\times ... \times [\frac{k_s}{2^{l_{i-1}}}\frac{k_s+1}{2^{l_{i-1}}} )$, we have from
 (\ref{3.56}),  (\ref{3.10}), (\ref{3.23}), Lemma 2 and (\ref{5.4}),  that
\begin{equation}\label{5.6}
      {\bf E} [ \cA(\dot{G}_i ) \; | \; \cF_{i-1}]  =O(n^{s} 2^{-   [n^{1/9}]} )=O(n^{-10s})   .
\end{equation}
Now let   $|\gamma_j| \leq 2^{(i+1)[n^{4/9}] -[n^{2/9}]}$   and $\theta_j^{(1,2)} \in  [\frac{k}{2^{l_i}},\frac{k+1}{2^{l_i}})$ then
\begin{equation}\nonumber
  | \sin(N_j \gamma_j \theta_j^{(1)}) - \sin(N_j \gamma_j \theta_j^{(2)}) |    
   = |N_j \gamma_j (\theta_j^{(1)} -\theta_j^{(2)}) \cos(N_j \gamma_j \theta_j^{(3)})| 
   \leq    2^{[n^{1/9}] -[n^{2/9}] +2} c_0 ,
\end{equation}
with  $\theta_j^{(3)} \in  [\frac{k}{2^{l_i}},\frac{k+1}{2^{l_i}})$, and 
\begin{equation}\nonumber
 \prod_{j=1}^s \sin(N_j \gamma_j \theta_j^{(1)}) =  \prod_{j=1}^s \sin(N_j \gamma_j \theta_j^{(2)})    
 +O(   2^{[n^{1/9}] -[n^{2/9}] } ).
\end{equation}
Therefore
\begin{equation}\label{5.9}
 \prod_{j=1}^s \sin(N_j \gamma_j \theta_j) =    
 {\bf E} \Big[ \prod_{j=1}^s \sin(N_j \gamma_j \theta) \; | \; \cF_{i}\Big] 
 +O(   2^{[n^{1/9}] -[n^{2/9}] } ).
\end{equation}
Taking into account  (\ref{3.56}),  (\ref{3.10}), (\ref{3.23})  and (\ref{5.6}), we get (\ref{5.2}). 
It is easy to see that
$$ 
|\cA(\dot{G}_i )^2 - \xi_i^2| \leq 
( 2|\cA(\dot{G}_i )| +  |\cA(\dot{G}_i ) - \xi_i|) | \cA(\dot{G}_i ) - \xi_i|, \quad and  \quad |\xi_i|^4  \leq 8|\cA(\dot{G}_i )|^4 +8  |\cA(\dot{G}_i ) - \xi_i|^4  .
$$
Applying (\ref{3.27b}), we obtain (\ref{5.21}). Hence, Lemma 18 is proved \qed
\\

 We shall use  the following variant of the {\it martingale central limit theorem} (see [Mo, p.~414]):
 
 Let $(\Omega, \cF,P)$ be a probability space and $\{(\zeta_{n,k}, F_{n,k}) \;| \; n=1,2,..., k=1,...,k_n\}$ be a martingale difference array with 
 ${\bf E} [\zeta_{n,k} | F_{n,k-1}] =0$ a.s. ($F_{n,0}$ is the trivial field).
 \\
 
{\bf Theorem  C.} {\it Let  $L(n,\epsilon) = \sum_{1 \leq k \leq k_n} {\bf E} [\zeta_{n,k}^2 \delta(|\zeta_{n,k}|>\epsilon)]$,
\begin{equation}\label{5.9a}
 \SS_{n} =\sum_{1 \leq k \leq i} \zeta_{n,k}, \qquad {\it and} \qquad \VV_{n}^2 =\sum_{1 \leq k \leq k_n} {\bf E} [\zeta_{n,k}^2 | F_{n,k-1}],
\end{equation}
\begin{equation}\label{5.90}
   \AA_n = {\bf E}[|\VV_{n}^2  -1|], \qquad  \WW_n= \int_0^1 L(n,\epsilon)  d \epsilon, \qquad {\rm and} \qquad \sum_{1 \leq k \leq k_n} 
   {\bf E} [\zeta_{n,k}^2] =1.
\end{equation}
Then }
\begin{equation}\label{5.92}
     \sup_t |P(\SS_{n} < t) -\Phi(t)|  \leq 7(\WW_n^{1/4}  + \AA_n^{1/3}).
\end{equation}
   Now we apply Theorem C to  the martingale difference array  (\ref{5.1}) with $F_{n,k}= \cF_k$, 
 $\zeta_{n,i} = \xi_i/\varrho$,  $\varrho= (\sum_{i \in [1,k_n]}  {\bf E} [\xi_i^2] )^{1/2}$, and  $k_n= [n^{5/9}]$.\\

 {\bf Lemma 19.} {\it Let 
\begin{equation}\label{5.94}
     \SS_n = \sum_{1 \leq i \leq k_n} \xi_i/\varrho.
\end{equation}    
Then }
\begin{equation}\nonumber
   \sup_t |P( \SS_n < t) -\Phi(t)|  = O(n^{-1/15}).
\end{equation}

{\bf Proof.}  By (\ref{5.1}),  $(\xi_i)_{i \geq 1}$ is the martingale difference sequence (and consequently orthogonal).
Using the triangle inequality,  Lemma 11, Lemma 13 and Lemma 18, we obtain
\begin{equation}  \nonumber
 \varrho^2 = \sum_{i \in [1,k_n]}  {\bf E} [\xi_i^2]  = {\bf E} \Big[\Big(\sum_{i \in [1,k_n]} \xi_i \Big)^2\Big] =
      {\bf E} \Big[\sum_{i \in [1,k_n]} \cA(\dot{G}_i )^2 \Big] +O(1)
\end{equation}  
\begin{equation}  \label{7.2}
                = {\bf E} [  (\cR (\btheta \cdot \bN  \cdot \KK^s+\bx, \Gamma) )^2  ] +O(n^{s-1-2/9})  \in n^{s-1} [w_3,w_4],
\end{equation}  
with some $w_4 > w_3 >0$.

Let $\dot{\cF}$  be a sub-$\sigma$-algebra of $\cF$. By Jensen's inequality, we get
\begin{equation}  \label{7.2a}
      {\bf E} [|\vartheta|^{\alpha}] \leq   ({\bf E} [|\vartheta|^{\beta}])^{\alpha /\beta} \quad {\rm and}  \quad
        {\bf E} [|\vartheta|^{\alpha} \; | \; \dot{\cF}] \leq   ({\bf E} [|\vartheta|^{\beta}\; | \; \dot{\cF}])^{\alpha / \beta}, \;\; {\rm with}  \;\; \beta >\alpha>0.
\end{equation}   
Consider $\WW_n$. 
 We derive from  (\ref{5.90}) that
\begin{equation}\nonumber
  \WW_n= \sum_{1 \leq i \leq k_n}   \int_0^1 \int  |\xi_i/\varrho|^{2} \delta(\xi_i/\varrho|> \epsilon) dP d \epsilon \leq \sum_{1 \leq i \leq k_n}   \int_0^1 \int_{ \{ |\xi_i/\varrho|> \epsilon \}} |\xi_i/\varrho|^{2} |\xi_i/(\varrho \epsilon)|^{24/25} dP d \epsilon
\end{equation} 
\begin{equation}\nonumber
  \leq    \sum_{1 \leq i \leq k_n}  \int |\xi_i/\varrho|^{74/25}  dP  \int_0^1\epsilon^{-24/25}  d \epsilon .
\end{equation}
Applying  (\ref{7.2}) and (\ref{7.2a}) with $\alpha=74/25$, $\beta=4$, we have
\begin{equation}\nonumber
      \WW_n \leq  25  \sum_{1 \leq i \leq k_n}  \Big( \int |\xi_i/\varrho|^{4}  dP \Big)^{37/50} =
  O\Big(n^{-(s-1)37/25} \sum_{1 \leq i \leq k_n} ({\bf E} [\cA^4(\dot{G}_i )])^{37/50}\Big) .
\end{equation}  
By Lemma 14, we get
\begin{equation}\nonumber
    \WW_n 
  =O\Big(n^{-(s-1)37/100} \sum_{1 \leq i \leq k_n} i^{2(s-2)\frac{37}{25}} n^{\frac{8}{9}(s-1)\frac{37}{50}} \Big)
\end{equation}  
\begin{equation}\label{7.4}
 =O(n^{-(s-1)\frac{37}{25} + \frac{5}{9}( 2(s-2)\frac{37}{50} +1 ) +\frac{8}{9}(s-1)\frac{37}{50}  }) =O(n^{-4/15})
   \;\; {\rm and} \;\;  \WW^{1/4}_n = O(n^{-1/15}). 
\end{equation} 
Next consider $\AA_n$. Let 
\begin{equation}\label{7.4b}
                   \UU_{n}^2 =\sum_{1 \leq k \leq k_n} (\xi_i/\varrho)^2.
\end{equation}   
Using (\ref{7.2}), (\ref{7.4b})  and  Lemma 18, we derive
\begin{equation}\nonumber
     {\bf E}[ |  \UU_n^2 -1|^2] =\varrho^{-4} {\bf E}\Big[\Big|\sum_{1 \leq k \leq k_n} (\xi_i^2 -  {\bf E}[\xi_i^2] )\Big|^2\Big]
\end{equation}   
\begin{equation}\nonumber
         \leq
   2  \varrho^{-4} {\bf E}\Big[\Big(\sum_{1 \leq k \leq k_n} (\cA^2(\dot{G}_i ) -  {\bf E}[\cA^2(\dot{G}_i )]) \Big)^2\Big] +O(n^{-5}).
\end{equation} 
By Lemma 17, we obtain 
\begin{equation}\label{7.10}
     {\bf E}[ |  \UU_n^2 -1|^2] =O(
    n^{-2(s-1) +2(s-1) -5/9} ) =O(n^{-5/9} ).
\end{equation} 
Let 
\begin{equation}\label{7.11}
   \varsigma_i = (\xi_i/\varrho)^2 - {\bf E} [(\xi_i/\varrho)^2 | \cF_{i-1}] \quad {\rm and } \quad  
   \VV_n^2 = \sum_{1 \leq i \leq k_n} {\bf E} [(\xi_i/\varrho)^2 | \cF_{i-1}].
\end{equation} 
By (\ref{5.1}), we see that
 $(\varsigma_i)_{i \geq 1}$ is the martingale difference array satisfying \\
   ${\bf E} [ \varsigma_i  |  \cF_{i-1}]=0, \; i=1,2,...$ .
From (\ref{7.4b}), we have
\begin{equation}\nonumber
 {\bf E}[ |\VV_n^2-\UU_n^2|^2] =  {\bf E}\Big[ \Big| \sum_{1 \leq i \leq k_n}  \varsigma_i  \Big|^2\Big] = \sum_{1 \leq i \leq k_n} {\bf E}[ \varsigma^2_i  ].
\end{equation}    
Using  (\ref{7.11}) and (\ref{7.2a}), we get
\begin{equation}\nonumber
     {\bf E}[ \varsigma^2_i  ]  \leq  2  \varrho^{-4}({\bf E}[ \xi^4_i]   +{\bf E} [({\bf E} [\xi^2_i | \cF_{i-1}])^2 ]) \leq 4\varrho^{-4}{\bf E}[ \xi^4_i  ].
\end{equation}  
By Lemma 14, Lemma 18 and (\ref{7.2}),    we have
\begin{equation}\label{7.00}
 {\bf E}[ |\VV_n^2-\UU_n^2|^2] = O\Big(n^{-2(s-1)}\sum_{1 \leq i \leq k_n}  {\bf E}  [\cA^4(\dot{G}_i )] +n^{-7s} \Big)
\end{equation} 
\begin{equation}\nonumber
= O\Big(n^{-2(s-1)}\sum_{1 \leq i \leq k_n}  i^{2(s-2)} n^{8/9(s-1)} \Big)=O(n^{-2(s-1) + 5/9(2s-3) + 8/9(s-1) })
  =O(n^{-5/9}).
\end{equation} 
 By (\ref{5.90}) and (\ref{7.2a}), we get
\begin{equation}\nonumber
  \AA^2_n= ({\bf E}[|\VV_n^2 -1|])^2 \leq {\bf E}[ |\VV_n^2 -1|^2]= {\bf E}[ |\VV_n^2-\UU_n^2  +  \UU_n^2 -1|^2] 
\end{equation}    
\begin{equation}\nonumber
 \leq 2  {\bf E}[ |\VV_n^2-\UU_n^2 |^2]
   + 2{\bf E}[|  \UU_n^2 -1|^2].
\end{equation} 
From (\ref{7.10}) and  (\ref{7.00}),    we derive
\begin{equation}\nonumber
  \AA^2_n=O(n^{-2/5} ),  \qquad {\rm and} \quad   \AA^{1/3}_n=O(n^{-1/15}).
\end{equation} 
Applying (\ref{7.4}) and Theorem C, we obtain the assertion of Lemma 19. \qed \\

\setcounter{equation}{0}
\renewcommand{\theequation}{2.8.\arabic{equation}}
{\bf 2.8. End of the proof of Theorem 1.}
\\
Let $\dot{\SS}_n =\cR( \boldsymbol{\theta} \cdot \bN \cdot \KK_s+\bx, \Gamma) / \dot{\varrho}$ and 
       $\dot{\varrho} =({\bf E}[\cR^2( \boldsymbol{\theta} \cdot \bN \cdot \KK_s+\bx, \Gamma) ])^{1/2}$.  
Using Lemma 11 and Lemma 18, we obtain
$$
     {\bf E}\Big[\Big( \cR( \boldsymbol{\theta} \cdot \bN \cdot \KK_s+\bx, \Gamma) -\sum_{1 \leq i \leq k_n} \xi_i\Big)^2\Big]  \leq 
             2 {\bf E}\Big[\Big( \cR( \boldsymbol{\theta} \cdot \bN \cdot \KK_s+\bx, \Gamma) - \sum_{1 \leq i \leq k_n}  \cA(\dot{G}_i )\Big)^2\Big] 
$$
\begin{equation}\label{8.4}
      +    2 {\bf E}\Big[\Big( \sum_{1 \leq i \leq k_n}( \cA(\dot{G}_i ) -  \xi_i)\Big)^2\Big] =O(n^{s-1-2/9}). 
\end{equation}
By (\ref{7.2}),  we get $ \dot{\varrho}^2 -\varrho^2  =O(n^{s-1-2/9}) $, $ \dot{\varrho}^2 \geq w_2 n^{s-1} $ for some $w_2 > 0$, and
\begin{equation}\label{8.2}
     \Big| \frac{1}{\varrho} - \frac{1}{\dot{\varrho}}\Big|
      = \frac{|\varrho- \dot{\varrho}|}{\varrho \dot{\varrho}} =  
       \frac{|\varrho^2- \dot{\varrho}^2|}{\varrho \dot{\varrho}|\varrho+ \dot{\varrho}|} =
     O(n^{  -3/2(s-1) -2/9}).
\end{equation}
Applying  (\ref{5.94}), (\ref{7.2}), (\ref{8.2}) and (\ref{8.4}), we derive
\begin{equation}\nonumber
   {\bf E}[(\SS_n- \dot{\SS}_n  )^2] \leq
                   2{\bf E} \Big[\Big( \sum_{1 \leq k \leq k_n} \xi_i -\cR( \boldsymbol{\theta} \cdot \bN \cdot \KK_s+\bx, \Gamma)  \Big)^2/\varrho^2\Big]
\end{equation}
\begin{equation}\nonumber
   + 2 ( 1/\varrho-1/\dot{\varrho})^2{\bf E}[ \cR^2( \boldsymbol{\theta} \cdot \bN \cdot \KK_s+\bx, \Gamma) ] =O(n^{-2/9}). 
\end{equation}
By Chebyshev's   inequality, we have
\begin{equation}\label{8.8}
  \qquad P(|\dot{\SS}_n  -\SS_n|  \geq   n^{-1/15}) = O(n^{-2/9+2/15}) =  O(n^{-1/15}).
\end{equation}
It is easy to see that 
$$
  \{ \dot{\SS}_n <t  \}   \subseteq \Big(\{ \SS_n <t +n^{-1/15} \} \cap \{ |\dot{\SS}_n -\SS_n| \leq   n^{-1/15} \} \Big) \cup  \{ |\dot{\SS}_n -\SS_n| \geq   n^{-1/15} \}
$$ 
and
$$
    \{ \SS_n <t -n^{-1/15} \}  \subseteq \Big(\{ \dot{\SS}_n <t  \} \cap \{ |\dot{\SS}_n -\SS_n| \leq   n^{-1/15} \} \Big) \cup  \{ |\dot{\SS}_n -\SS_n| \geq   n^{-1/15} \}.
$$
Hence
$$
  P(\{ \SS_n <t -n^{-1/15} \} ) -P(\{ |\dot{\SS}_n -\SS_n| \geq   n^{-1/15} \})       \leq  P(\{ \dot{\SS}_n <t  \})  
$$
\begin{equation}\label{8.10}
            \leq  P(\{ \SS_n <t +n^{-1/15} \} )      +P(\{ |\dot{\SS}_n -\SS_n| \geq   n^{-1/15} \}).
\end{equation}
We note for $|u| \leq n^{-1/15}$ that
\begin{equation}\nonumber
  |\Phi(t +u) -\Phi(t) |  <  \frac{1}{\sqrt{2 \pi}} \int_{t- |u|}^{t+|u|}  e^{-u^2/2}du
                \leq \frac{1}{\sqrt{2 \pi}} \int_{t- n^{-1/15} }^{t+n^{-1/15}} du = \frac{2}{\sqrt{2 \pi}} n^{-1/15}.
\end{equation}
Using Lemma 19,  we get  
\begin{equation}\nonumber
     \sup_t | P(\{ \SS_n <t +u \} )  -\Phi(t)|  \leq  \sup_t (| P(\{ \SS_n <t+u \} )  -\Phi(t +u)| 
\end{equation}
\begin{equation}\nonumber
     +|\Phi(t +u)-\Phi(t) |) =O(n^{-1/15}  ),  \quad  |u| \leq n^{-1/15}.
\end{equation}
By (\ref{8.10}) and (\ref{8.8}), we derive
\begin{equation}\nonumber
   \sup_t |P( \dot{\SS}_n < t) -\Phi(t)|  =O(n^{-1/15}).
\end{equation}

Bearing in mind that throughout the paper $O$-constants does not depend on $\bx$,
we obtain the assertion of Theorem 1. \qed
\\

\setcounter{equation}{0}
\renewcommand{\theequation}{2.9.\arabic{equation}}

{\bf 2.9. Sketch of the proof of Theorem 2.} 
We use notations from \S 1.3. Let \\  $I_0=[0,y_1) \times \cdots \times [0,y_{s-1}) $, $I_1=[-y_1/2, y_1/2) \times \cdots \times [-y_{s-1}/2,y_{s-1}/2) $, $I_2 = [-1/2,1/2)^{s-1 }$, \\
  $I_3=[-y_{s} N \det \Gamma/2,     y_{s} N \det \Gamma/2) $,   $I_4= [-z_s(\bx, [y_s N] )/2, z_s(\bx, [y_s N] )/2)$,   \\  $\bu_1=(y_1/2,...,y_{s-1}/2, z_s(\bx, [y_s N] )/2) -\dot{\bx}$, $\bu_2=(1/2,...,1/2, z_s(\bx, [y_s N] )/2) -\dot{\bx}$  with $\bx=(x_1,...,x_{s-1})$, 
   $\dot{\bx}=(x_1,...,x_{s-1},0)$.
   
By (\ref{1.0}), (\ref{N1}) and (\ref{N2}), we obtain  
$$
  \Delta(I_0, (\cT^l(\bx))_{l=0}^{[y_s N]-1})    = 
 \cN(  I_1  \times I_3 +\bu_1, \Gamma)
 - y_1 \cdots y_{s-1} \cN(  I_3  \times I_2 +\bu_2, \Gamma).
$$
Let $a= z_s(\bx, [y_s N] )$, $b =y_{s} N \det \Gamma$, and let
\begin{equation} \nonumber
   (\kappa,I) =   \begin{cases} 
    (1, [-a/2,-b/2) \cup [b/2,a/2),  & \; {\rm if}  \;  a>b,\\
    ( -1, [-b/2,-a/2) \cup [a/2,b/2 , &{\rm otherwise}.
  \end{cases}
\end{equation}
By (\ref{N1}) and (\ref{N2}),  we get
$$
  \Delta(I_0, (\cT^l(\bx))_{l=0}^{[y_s N]-1})    =  \dot{R_1} + \kappa \ddot{R_1} - y_1 y_2\cdots y_{s-1} (\dot{R_2} + \kappa \ddot{R_2}),
$$
with
$$
   \dot{R_i} =  \cR( I_i \times I_3       +\bu_k, \Gamma), \quad  {\rm and}  \quad \ddot{R_i} =  \cR( I_i \times I_4      +\bu_k, \Gamma), \quad i=1,2.
$$
It is easy to verify (see also [Le2, p. 86]) that
$$
     \ddot{R_i} =O((\ln(n))^{s-1}), \quad i=1,2.
$$
Thus $\dot{R_1}- y_1 y_2\cdots y_{s-1}\dot{R_2}$ is the essential part of $\Delta(I_0, (\cT^l(\bx))_{l=0}^{[y_s N]-1})  $.
Repeating the proofs of $\S 2.4$,  we have the   upper bound of the variance of $\dot{R_1}- y_1 \cdots y_{s-1}\dot{R_2}$.  Using Roth's inequality (\ref{1.6}),
 we get the lower bound of the variance $\dot{R_1}- y_1 \cdots y_{s-1}\dot{R_2}$. Next repeating the proofs of $\S 2.5-\S2.8$,  we obtain the assertion of Theorem 2. \hfill \qed  
\\

{\bf Bibliography.}
  
  [BW] Baker, A., Wustholz, G., Logarithmic Forms and Diophantine Geometry.  Cambridge University Press, Cambridge, 2007. 

 [Be1] Beck, J.,    Randomness of $n\sqrt 2\bmod 1$ and a Ramsey property of the hyperbola. Sets, graphs and numbers (Budapest, 1991), 23--66, Colloq. Math. Soc. Janos Bolyai, 60, North-Holland, Amsterdam, 1992. 

  [Be2]  Beck, J., Randomness in lattice point problems, Discrete Math. 229 (2001), no. 1-3, 29-55.

   [Be3]     Beck, J. Randomness of the square root of 2 and the giant leap, Parts 1,2, Period. Math. Hungar. 60 (2010), no. 2, 137-242; 
 62 (2011), no. 2, 127-246

[BC] Beck, J., Chen, W.W.L.,
 Irregularities of Distribution,  
 Cambridge Univ. Press,  Cambridge 1987. 

[BS]  Borevich, A. I., Shafarevich, I. R., Number Theory, Academic Press, New York, 1966. 

[By] Bykovski\u\i, V.A., On the right order of error of optimal cubature formulas in the spaces with dominating derivation and $L^2$ discrepancy of nets, Dalnevost. Science Center of the USSR Acad. of Sciences,
  Vladivostok, 1985 (in Russian).

[DrTi]  Drmota, M.,  Tichy, R., Sequences, Discrepancies and Applications, Lecture Notes in Mathematics
1651, 1997.

  [ESS]  Evertse, J.H., Schlickewei, H.P., Schmidt, W.M., Linear equations in variables which lie in a multiplicative group, Ann. of Math. (2) 155 (2002), no. 3, 807-836.

[Fr] Frolov, K. K., Upper bound of the discrepancy in metric $L\sb{p}$, $2\leq p<\infty $. (Russian) Dokl. Akad. Nauk
 SSSR 252 (1980), no. 4, 805-807.
 [English translation: Soviet Math. Dokl. 21 (1980), no. 3, 840-842 (1981).]

 [GL]  Gruber, P.M,  Lekkerkerker, C.G., Geometry of Numbers, North-Holland, New-York,   1987.
 
 [HuRu] Hughes, C. P., Rudnick, Z., On the distribution of lattice points in thin annuli, Int. Math. Res. Not.   13 (2004), 637-658.

  [Le1]  Levin, M. B., The multidimensional generalization of J.Beck 'Randomness of $n \sqrt 2 \;
  \mod 1$ ...'  and a.s. invariance principle for $\ZZ^d$-actions of toral 
 automorphisms, Abstracts of Annual Meeting of the Israel Mathematical Union,(2002), \\ http://imu.org.il/Meetings/IMUmeeting2002/ergodic.txt.  
 [Le2] Levin, M.B., On low discrepancy sequences and low
discrepancy ergodic transformations of the multidimensional unit
cube,    Israel J.  Math. {\bf 178} (2010),  61-106.

  [Le3]  Levin, M.B.,  Central Limit Theorem for $\ZZ_{+}^d$-actions  by toral endomorphisms,  Electronic Journal of Probability 18 (2013), no. 35, 42~pp.

 [LeMe]  Levin, M.B.,  Merzbach, E., Central limit theorems for the ergodic
adding machine, Israel J.  Math. 134 (2003),  61-92.  
  
[Mo] Mori, T. On the rate of convergence in the martingale central limit theorem. Studia Sci. Math. Hungar. 12 (1977), no. 3-4, 413–417.  
  
 [Skr] Skriganov, M.M.,  Construction of uniform distributions in terms of
geometry of numbers,
 Algebra i Analiz {\bf6}, no.~3 (1994), 200-230; 
Reprinted in St.~Petersburg Math. J., {\bf6}, no.~3 (1995), 635-664.

   [SW]  Stein, E., Weiss, G.,  Introduction to Fourier Analysis on Euclidean Spaces, Princeton University Press, New-York, 1971.
\\
\\
{\bf Address}: Department of Mathematics,
Bar-Ilan University, Ramat-Gan, 52900, Israel \\
{\bf E-mail}: mlevin@math.biu.ac.il\\
\end{document}